\documentclass[11pt,a4paper]{article}
\usepackage{geometry,setspace}
\usepackage{euler}
\usepackage{pict2e,picture}

\makeatletter
\DeclareRobustCommand{\bbDelta}{{\mathpalette\bb@Delta\relax}}
\newcommand{\bb@Delta}[2]{%
  \begingroup
  \sbox\z@{$\m@th#1\Delta$}%
  \dimendef\Dht=6 \dimendef\Dwd=8
  \setlength{\Dwd}{\wd\z@}%
  \setlength{\Dht}{\ht\z@}%
  \begin{picture}(\Dwd,\Dht)
  \put(0,0){$\m@th#1\Delta$}
  \put(.42\Dwd,.7\Dht){\line(10,-26){.25\Dht}}
  \end{picture}%
  \endgroup
}
\makeatother

\usepackage[utf8]{inputenc}
\usepackage{amsmath}
\usepackage{amssymb}
\usepackage{mathtools}
\usepackage{amsfonts}
\usepackage{graphicx}
\usepackage{xcolor}
\usepackage{hyperref}
\usepackage{float} 
\usepackage[ruled,vlined]{algorithm2e}
\newcommand{\cL}{\mathcal{L}}
\newcommand{\cQ}{\mathcal{Q}}

\newcommand{\cH}{\mathcal{H}}
\newcommand{\cS}{\mathcal{S}}

\newcommand{\X}{\mathbb{X}}
\newcommand{\R}{\mathbb{R}}
\newcommand{\Z}{\mathbb{Z}}

\newcommand{\hsigma}{\tilde{\sigma}}

\newcommand{\tphi}{\tilde{\varphi}}
\newcommand{\tla}{\tilde{\lambda}}
\newcommand{\tLa}{\tilde{\Lambda}}

\newcommand{\tK}{\tilde{K}}
\newcommand{\tA}{\tilde{A}}
\newcommand{\tB}{\tilde{B}}
\newcommand{\tL}{\tilde{\mathcal{L}}}
\newcommand{\tLast}{\tL^{\ast}}
\newcommand{\tf}{\tilde{f}}
\newcommand{\trho}{\tilde{\rho}}
\newcommand{\tmu}{\tilde{\mu}}
\newcommand{\dt}{\mathrm{d}t}
\newcommand{\dW}{\mathrm{d}W}

\newcommand{\EE}{\mathbb{E}} % expectation value
\newcommand{\inv}{^{-1}}

\newcommand{\tn}{{\theta_n}}
\usepackage{pgfplots}
\pgfplotsset{compat=1.18}
\usepgfplotslibrary{groupplots}
\usepackage[affil-it]{authblk} 

\title{Effective Dynamics and Transition Pathways from Koopman-Inspired Neural Learning of Collective Variables}

\author[1,2]{Alexander Sikorski}
\author[2,1]{Luca Donati}
\author[1]{Marcus Weber}
\author[1,2,*]{Christof Schütte}

\affil[1]{Zuse Institute Berlin (ZIB), Takustr. 7, 14195 Berlin, Germany}
\affil[2]{Freie Universität Berlin, Department of Mathematics and Computer Science, Arnimallee 6, 14195 Berlin, Germany}
\affil[*]{schuette@zib.de}

\date{}

\begin{document}
\maketitle

\begin{abstract}
The ISOKANN (Invariant Subspaces of Koopman Operators Learned by Artificial Neural Networks) framework provides a data-driven route to extract collective variables (CVs) and effective dynamics from complex molecular systems. In this work, we integrate the theoretical foundation of Koopman operators with Krylov-like subspace algorithms, and reduced dynamical modeling to build a coherent picture of how to describe metastable transitions in high-dimensional systems based on CVs. Starting from the identification of CVs based on dominant invariant subspaces, we derive the corresponding effective dynamics on the latent space and connect these to transition rates and times, committor functions, and transition pathways. The combination of Koopman-based learning and reduced-dimensional effective dynamics yields a principled framework for computing transition rates and pathways from simulation data. Numerical experiments on one-, two-, and three-dimensional benchmark potentials illustrate the ability of ISOKANN to reconstruct the coarse-grained kinetics and reproduce transition times across enthalpic and entropic barriers.
\end{abstract}

% \section*{Alternative Title}
% \begin{enumerate}
%     \item \textbf{Learning Collective Variables and Reduced Dynamics via ISOKANN: Theory, Algorithms, and Applications}
% \end{enumerate}

%\newpage
%\tableofcontents
\newpage

\section{Introduction}

The exploration of high-dimensional stochastic dynamical systems, such as those arising in molecular dynamics, is often hindered by the complexity and dimensionality of the underlying processes. Effective low-dimensional descriptions are therefore indispensable for understanding, simulating, and interpreting the dominant mechanisms that govern such systems. The notion of \emph{collective variables} (CVs) provides a natural framework for this reduction: by identifying a few appropriately chosen variables that capture the essential slow degrees of freedom, one may represent the dynamics of the full system through an effective, coarse-grained model.

Traditional approaches to coarse-graining rely on physical intuition or prior knowledge of relevant coordinates. In recent years, however, \emph{algorithmic identification of CVs} has emerged as a promising alternative, exploiting advances in data-driven and machine-learning techniques to discover optimal reaction coordinates directly from simulation data. Examples include time-lagged independent component analysis (TICA) \cite{PerezHernandez2013TICA}, variational approaches to Markov processes (VAMP, VAMPnets)~\cite{Mardt2018VAMPnets}, approaches based on normalizing flows \cite{WuNoe2024ReactionCoordinateFlows} and dominant subspace learning methods such as ISOKANN \cite{ISOKANN}, which minimize the deviation between effective and full dynamics in the sense of stochastic model reduction~\cite{ZhangSchuette2025,NueskeKoltaiBoninsegnaClementi2021}. These approaches extend earlier theoretical developments on optimal reaction coordinates and exact coarse-graining without time-scale separation~\cite{LuVandenEijnden2014ExactCoarseGraining}, as well as mathematical analyses of pathwise error bounds and effective dynamics~\cite{ZhangHartmannSchuette2016,LegollLelievreOlla2017Pathwise,LelievreZhang2019Pathwise}.

At the same time, rigorous frameworks for connecting microscopic and effective dynamics have clarified under what conditions a reduced dynamics provides a faithful description of the original stochastic process~\cite{Acta}. Yet, despite these advances and some algorithmic progress \cite{Koumoutsakos2022}, a unified understanding of how the \emph{data-driven identification} of CVs interacts with the \emph{mathematical and algorithmic characterization} of their effective dynamics remains incomplete.

The aim of this work is to demonstrate how these two aspects---the algorithmic determination of collective variables and the analysis of their effective dynamics---can be pursued jointly and consistently. We illustrate this interplay using simple yet representative examples, showing explicitly how the main transition rates and pathways of the effective dynamics reconstructed from the learned CVs compare to the rates and pathways of the underlying full dynamics. Through these case studies, we provide evidence that data-driven discovery of CVs, when guided by principles from stochastic reduction theory, can lead to quantitatively reliable and physically interpretable coarse-grained models.

In the remainder of this article, we first review the theoretical foundations of collective variables (CVs) including their relation to the spectral properties of the Koopman operator as well as the fundamentals of the associated effective dynamics (Section~\ref{sec:CV}). Building on this, Section 3 introduces the ISOKANN framework, detailing both its mathematical underpinnings and its practical algorithmic realization through neural networks and the inner simplex transformation. Section~\ref{sec:transitions} then develops the formulation of effective transition kinetics in both the full and the reduced dynamics, including committor functions, reactive fluxes, and transition rates. Section~\ref{sec:numerical_experiments} presents a sequence of numerical experiments in one-, two-, and three-dimensional systems that illustrate how ISOKANN recovers dominant eigenmodes, coarse-grained potentials, diffusion tensors, and transition pathways. Finally, Section~\ref{sec:conclusion} summarizes the conceptual insights and outlines directions for future research on data-driven identification of CVs and effective stochastic models.

%%%%%%%%%%%%%%%%%%%%%%%%%%%%%%%%%%%%
\section{Collective Variables (CVs) in Molecular Dynamics}\label{sec:CV}

In this section, the transfer operator theory for CVs is shortly introduced and connected with the relation between slow modes of the full dynamics with the dominant eigenmodes of the Koopman operator and the effective dynamics that is induced by a CV on its latent space.  

\subsection{Koopman Operators, Dominant Eigenvalues and Slow Modes}

For a stochastic process $(x_t)_{t\in\R}$ depending on time $t$ with states $x_t\in \X\subseteq \R^d$ (e.g., MD), the Koopman operator $K_\tau$ is a transfer operator of the process. It acts on functions $f$ as
\[
\label{eq:defkoopman}
(K_\tau f)(x) = \EE[f(x_{t+\tau}) \mid x_t = x].
\]
We assume that the process $x_t$ is ergodic and admits a unique invariant measure $\mu$ with nonnegative density $\mu(x)$ such that $\mu_A=\int_A \mu(x)dx$ for each measurable set $A\subset \X$.  
We consider $K_\tau$ as an operator acting on the weighted Hilbert space $\cH=L^2_\mu$, the linear function space of square integrable functions, equipped with scalar product $\langle f,g\rangle_\mu =\int f(x)g(x)\mu(x)dx$.  
That is, $K_\tau$ is a linear (but infinite-dimensional) operator whose eigenfunctions $\varphi\in \cH$ satisfy:
\[
K_\tau \varphi = \Lambda\varphi
\]
where $\Lambda_i \in \mathbb{C}$. 

Subsequently, we assume that the process $x_t$ is reversible and $\mu$ is absolutely continuous. This is an assumption made for simplicity, the general case is described in detail in \cite{Acta}. Under this assumption, $K_\tau$ is self-adjoint and its eigenvalues are real-valued and discrete. Let $\Lambda_0, \Lambda_1, \ldots$ denote the Koopman eigenvalues in descending order and let $\varphi_i$ denote the eigenfunction of $\Lambda_i$ such that $\langle\varphi_i,\varphi_i\rangle_\mu=\delta_{ij}$.  We always have $\Lambda_0=1$ with $\varphi_0 = 1$, and that
any function $f\in\cH$ can be written as
\[
f = \sum_{i=0}^\infty \langle f, \varphi_i\rangle_\mu \varphi_i.
\]
Then, the Koopman operator $K_\tau$ admits eigenfunctions whose spectral properties reflect slow transitions between macrostates:
\begin{equation}\label{eq:eigenexpansion}
K_\tau^n f=\sum_{i=0}^\infty \Lambda_i^n\langle f, \varphi_i\rangle_\mu \varphi_i,
\end{equation}
such that information coded in the \emph{slow modes} or \emph{dominant macrostates} $\varphi_i$, the eigenfunctions of the dominant first eigenvalues $\Lambda$ close to $\Lambda_0=1$, decays slowest and information of higher eigenfunctions with eigenvalues $\Lambda\ll 1$ decays exponentially fast. 

If the eigenvalues of $K_\tau$ admit a spectral gap, that is, 
\[
1=\Lambda_0>\Lambda_1\ge \ldots\Lambda_{m}\gg \Lambda_{m+1}
\]
then the $m+1$ dominant eigenvalues $\Lambda_0,\ldots,\Lambda_{m}$ and the associated eigenfunctions are called the \emph{slow modes} and the process is called \emph{metastable}. In this case, the process exhibits $m+1$ distinct metastable sets $A$, each having the property that exit from the set is a rare event, that is, exit from the set takes very long on average and happens on timescales $\tau/|\log(\Lambda)|$, cf.~\cite{Acta}.

For example, the Markov process $(X_t)$ given by the diffusive molecular dynamics process 
\begin{equation} \label{eq:full}
d X_t = - \nabla V(X_t) \dt +\sigma \dW_t,
\end{equation}
with smooth potential $V$ that is growing to infinity fast enough, constant noise intensity $\sigma>0$ and $d$-dimensional Wiener noise process $W_t$, is reversible with absolutely continuous invariant density $\mu(x)\propto \exp(-\beta V(x))$, where $\beta=2/\sigma^2$. The evolution of expectation values, $f(x,t)=\mathbb{E}[f(X_{t+\tau})| X_0=x]$, under this dynamics is described by 
\[
\partial_t f(x,t) =\cL f(x,t),
\]
where the generator of the process is given by 
\begin{equation} \label{eq:full_generator}
\cL f = \frac{\sigma^2}{2}\Delta f - \nabla V\cdot \nabla f,
\end{equation}
such that the Koopman operator satisfies
\[
K_\tau =\exp(\tau\cL).
\]
Since the dynamics is reversible with invariant density $\mu$, $\cL$ is self-adjoint in $\cH$, i.e., it is identical to its adjoint $\cL^\ast$ and the Koopman operator coincides with the evolution operator of the operator $\cL^\ast$.
That is, in $\cH$, the evolution equation for probability densities under the dynamics, reads
\[ 
\partial_t\rho=\cL^\ast \rho=\cL\rho,
\] 
where $\rho$ is defined as the probability density with respect to the invariant density $\mu$ with $\rho \rightarrow 1$ as $t\rightarrow \infty$, indicating convergence of the system toward equilibrium.
Note that $\rho$ should not be confused with the probability density with respect to the Lebesgue measure solution of the usual Fokker-Planck equation in $L^2(dx)$.

$K_\tau$ and $\cL$ have the same eigenfunctions and the eigenvalues $\lambda$ of $\cL$ are related to the Koopman eigenvalues $\Lambda$ by $\Lambda=\exp(\tau\lambda)$, showing that the lowest eigenvalues $\lambda\approx 0$ of $\cL$ correspond to the dominant eigenvalues $\Lambda\approx 1$ of $K_\tau$. That is, if there are only two dominant eigenvalues $\lambda_0=0$ and $\lambda_1\approx 0$ with two distinct metastable sets, the expected exit from these sets happens on timescale $1/|\lambda_1|$, see \cite{Bovier1,Acta}. 

The \emph{transfer operator theory} was first introduced for molecular dynamics in \cite{SFHD99}, and has led to several hundreds articles about specific theoretical approaches \cite{SS14,Acta}, data-based approximation as discussed in \cite{KNKWKSN18}, Markov State Models (MSMs), as reviewed in \cite{bowman2013}, variational approaches utilizing neural nets as in \cite{Mardt2018VAMPnets}, and many other directions.

\subsection{Collective Variables}
The notion of \emph{collective variables} (CV), often also called \emph{reaction coordinates}, is used in different versions in the literature, and considerable theory has been developed recently. The basic concept is that the collective variable is an abstract, low-dimensional set of coordinates which represents progress along transition pathways between an initial state and a target state that mostly belong to different metastable sets. In fact, there are two slightly different ways to understand the last sentence: (I) A CV is understood as a classical reaction coordinate, i.e., it describes the progress of a specific transition between two metastable sets and therefore is crafted for this process in particular. In this case, the committor function, as introduced in Sec.~\ref{sec:transitions} below, has been demonstrated to be the optimal reaction coordinate \cite{LuVandenEijnden2014ExactCoarseGraining,E2005Transition}. (II) The CV allows for a description of many such transition processes (the ``important'' ones), is ``global'' for the system under investigation, and provides a parametrization of the \emph{transition manifold} of the system as described in \cite{Bittracher2018}. We will herein adopt (II).   

In general, a CV is a nonlinear, smooth map $ \xi \colon \X \to \R^m $ that reduces the dimension $d=\mathrm{dim}(\X)$ to a significantly smaller dimension $m\ll d$ with the main additional requirement that the ``projection'' of full-dimensional molecular dynamics to the CV allows for a ``good'' reproduction of the long-term dynamical behavior of the system under investigation. We call $\Z=\xi(\X)$ the \emph{latent space} of the CV.

Let $ \xi \colon \X \to \R^m $ be a smooth ($C^1$) function, where $ m \le d $, and let $ \mathbb{L}_z = \{x \in \X \mid \xi(x) = z\} $ be the $ z $-level set of $ \xi $, that is, $\mathbb{L}_z=\{x\in\X:\,\xi(x)=z\}$. The coordinate projection $\Pi_\xi$ averages a given function along the level sets of a coordinate function~$\xi$: 
\begin{eqnarray}
    \Pi_\xi f(x)  & = &  \mathbb{E}_\mu(f(x')|\xi(x')=\xi(x)).\label{eq:Pxi2}
\end{eqnarray}
That is, $\Pi_\xi f$ is constant along the level sets $\mathbb{L}_z$, and is defined by the average of $f$ along these level sets. Therefore, $\Pi_\xi f$ is often also understood as a function $z$ on the latent space, $z\in\xi(\X)$, in the sense of $\tilde{f}(z)=(\Pi_\xi f)(\xi^{-1}(z))$, with
\[
\Pi_\xi f(x)=\int_{\mathbb{L}_{\xi(x)}} f(x)\mu_{\xi(x)}(x)dx,
\]
where~$ \mu_z=\mu_{\xi(x)} $ is a probability measure on $ \mathbb{L}_z $, the marginal of $\mu$ conditional to $ \xi(x) = z$, such that
\[
\mathbb{E}_\mu(f)=\int_{\X} f(x)\mu(dx)= \int_{\Z} \int_{\mathbb{L}_z} f(x)\mu_z(dx)\; \tilde{\mu}(z)dz=\int_{\Z} \mathbb{E}_\mu(f(x)|\xi(x)=z)\; \tilde{\mu}(z)dz,
\]
where $\tilde{\mu}(z)$ is the measure induced by $\mu$ on the latent space $\Z=\xi(\X)$.
Because of this, one often simplifies notation and writes $\Pi_\xi f(z)$ as function in latent space $\mathbb{Z}$ which represents $\xi(\X)$ via the coordinate representation $z=\xi(x)$.   For the canonical measure $\mu(x) = Z^{-1}\exp(-\beta V(x))$, a particular case is given by the function $f=1$ with $\Pi_\xi 1(z)=\int_{\mathbb{L}_z}1(x)\mu_z(x) d x$, which allows to denote the \emph{free energy} function associated with the CV,
\begin{equation}\label{free energy}
F_\xi(z)=-\frac{1}{\beta}\log \Pi_\xi 1(z).
\end{equation}
The full process $(X_t)$ defines a dynamics in $m$-dimensional \emph{latent space} spanned by $\xi$ by 
\[
z_t=\xi(X_t),
\]
and an application of Ito's lemma yields, based on (\ref{eq:full}),
\begin{equation}\label{eq:projected}
   d z_t  =  (\cL \xi)(X_t) \dt +\sigma \nabla\xi(X_t)\cdot \dW_t,
\end{equation}
which is \emph{not closed} since the right hand side not only depends on the value
of the collective variable, $z_t$, but also on the state of the original full dynamics, $X_t$. For circumventing this problem,
there are mainly two options: (A) use elaborated closing relations for $X_t$ by introducing memory terms that aim to approximate the correlations missing by terms only depending on $z_t$, or (B) bring (\ref{eq:projected}) into a closed Markovian form by replacing the coefficients in (\ref{eq:projected}) with their projections via $\Pi_\xi$ to the latent space. While (A) is used in the literature mainly via so-called generalized Langevin equations \cite{Netz2018,Ayaz2021NonMarkovian} or learning the effective dynamics via Long short-term memory (LSTM) \cite{Koumoutsakos2022}, we will herein follow approach (B) which has attracted a lot of attention in the literature recently \cite{Legoll2010,ZhangHartmannSchuette2016,ZhangSchuette2017,LegollLelievreOlla2017Pathwise,LegollLelievreSharma2019Effective,NueskeKoltaiBoninsegnaClementi2021,ZhangSchuette2025}.
The specific approach utilized herein stems from \cite{Bittracher2018} and is based on the following key idea: Find CVs $\xi$ for which the Markovian closed form, called the \emph{effective dynamics} induced by $\xi$, exhibits (approximately) the same slow modes than the original dynamics. 

\subsection{Effective Dynamics and Inherited Slow Modes}
The \emph{effective dynamics} given by the CV $\xi$ is the dynamics in latent space (dimension $m$) induced by the full dynamics  (\ref{eq:full}). That is, the full dynamics is projected onto the latent space via the projection operator $\Pi_\xi$. In \cite{ZhangSchuette2025}, the effective dynamics is discussed in the general case (not only for diffusion processes). For simplicity of presentation, we here restrict ourselves to the case that the dynamics is given by (\ref{eq:full}).
According to \cite{ZhangHartmannSchuette2016,ZhangSchuette2017}, in this case, the effective dynamics is then given by 
\begin{align}
  dz_t = \widetilde{b}(z_t) \,\dt + \sigma\cdot
  A(z_t) \,\dW_t\,,
  \label{eff-dynamics}
\end{align}
where $z_t$ denote the state on the latent space at time $t$.  The coefficients $\widetilde{b} : \mathbb{R}^{m}\rightarrow
\mathbb{R}^{m}$, $A:
\mathbb{R}^{m} \rightarrow \mathbb{R}^{m \times m}$ live on the latent space $\Z$ (no longer full state space), together with the process $W_t$, with 
\begin{align}
  \begin{split}
    \widetilde{b}_l(z) =& \Pi_\xi(\mathcal{L}\xi_l)(z), \\
    (AA^\top)_{lk}(z) &= \Pi_\xi\Big(
    \sum_{i,j=1}^{d}  \frac{\partial \xi_l}{\partial x_i} \frac{\partial \xi_k}{\partial x_j}
    \Big)(z),
\end{split}
\label{effective-coeff}
\end{align}
for $z \in \R^m$, so we can assume that $A=A(z)$ is an $m \times m$ square matrix. We also have the following alternative formula
\begin{align}
  \begin{split}
    \widetilde{b}_l(z) =& \lim_{s \rightarrow 0+}
    \EE\Big(\frac{\xi_l(X_s)-z_l}{s}~\Big|~X_0 \sim \mu_z\Big),\\
    (AA^\top)_{lk}(z) =&
    \lim_{s \rightarrow 0+} \EE
  \Big( \frac{(\xi_l(X_s) - z_l)(\xi_k(X_s) - z_k)^\top}{s}~\Big|~X_0 \sim
  \mu_z\Big),
\end{split}
\label{effective-coeff-alternative}
\end{align}
where the conditional expectations are with respect to the ensemble of trajectories of the full
dynamics (\ref{eq:full}) with the initial distribution $\mu_z$ on
$\mathbb{L}_z$.

The generator $\tL$ of (\ref{eff-dynamics}) results from the generator $\cL$ of the full process being projected to latent space, 
\[
\Pi_\xi \cL f(x) = (\tL \tf)(\xi(x)),
\]
with $f(x)=\tf(\xi(x))$, which is given by
%
% \begin{equation}\label{eff-gen}
% \left(\tL \tf\right)(z) =\widetilde{b} \cdot \nabla_z \tf 
% + 
% \frac{\sigma^2}{2}\text{Tr}\left(
% AA^\top(z)\nabla_z^2 \tf \right),
% \end{equation}

\begin{equation}\label{eff-gen}
\tL = \sum_{l=1}^m\widetilde{b}_l \frac{\partial}{\partial z_l}+\frac{1}{2}\sum_{k,l=1}^m (\hsigma\hsigma^\top)_{kl}(z) \frac{\partial^2}{\partial z_k\partial z_l},
\end{equation}
where $\hsigma\hsigma^\top=\sigma^2 AA^\top$.
The \emph{effective Koopman operator} then is defined by
\[
\tK_\tau=\exp(\tau \tL).
\]
Like the effective generator $\tL$, $\tK_\tau$ acts on functions on the latent space and is the latent space variant of the projected Koopman operator $\Pi_\xi K_\tau \Pi_\xi$. If $K_\tau$ is self-adjoint then $\tK_\tau$ is too (on the correctly weighted Hilbert space $\tilde{\mathcal{H}}$ on $\xi(\X)$).

Concerning the eigenvalues of $\tK_\tau$ we have the 
following statement \cite{ZhangSchuette2025} (Thm. 4.1): For each eigenvalue $\Lambda$ and eigenfunction $\varphi$ of $K_\tau$ we have that:
\begin{itemize}
 \item[(1)] For each eigenvalue $\tLa$ of $\tK_\tau$ with eigenfunction $\tphi$ we have the estimate
 \[
 \Lambda-\tLa\le \mathcal{E}(\varphi-\tphi\circ\xi)-(1-\Lambda)\|\varphi-\tphi\circ\xi\|_\mu,
 \]
with $\mathcal{E}(f)=\langle(\text{Id}-K_\tau)f,f\rangle_\mu$ where $\circ$ denote concatenation, i.e., $(\tphi\circ\xi)(x)=\tphi(\xi(x))$.
 \item[(2)] If there is a function $\tphi$ on the latent space such that $\varphi(x)=\tphi(\xi(x))$, then $\Lambda$ is also an eigenvalue of $\tK_\tau$ with eigenfunction $\tphi$.
 \item[(3)] When ordering the spectrum of $K_\tau$ and $\tK_\tau$ in decreasing order, we have
 \begin{equation}\label{order}
 \Lambda_i\ge\tLa_i.
 \end{equation}
\end{itemize}
Consequently, whenever we want that the slow modes of $\tK_\tau$ are similar to the ones of $K_\tau$, we have to look for CVs $\xi$ for which the dominant eigenfunction $\varphi$ are almost constant on the level sets $\mathbb{L}_z$, because then we will find a $\tphi$ such that $\varphi(x)\approx \tphi(\xi(x))$. 

Given this criterion, the obvious choice for a CV of dimension $m$ would be given by the first nontrivial dominant eigenfunctions, 
\[
\xi(x)=(\varphi_1(x),\ldots,\varphi_m(x)),
\]
since, then, the slow modes of the full dynamics, i.e., dominant eigenvalues of $K_\tau$, were inherited by the effective dynamics, i.e., the dominant eigenvalues of $\tK_\tau$. However, algorithmically, this choice is problematic since the eigenfunctions of $K_\tau$ are unknown and computation is difficult. Nevertheless, there is a growing variety of numerical approaches using neural networks and deep learning for identifying the dominant eigenfunctions, e.g., see \cite{Mardt2018VAMPnets,ZhangLiSchuette2022EigenvaluePDEs}, or the subspace spanned by them \cite{ISOKANN}.

%%%%%%%%%%%%%%%%%%%%%%%%%%%%%%%%%%%%%%%%%%%%%%%%%%%%
\section{Identification of CVs via ISOKANN}

\emph{ISOKANN} stands for \textit{Invariant Subspaces of Koopman Operators Learned by Artificial Neural Networks}. It is a data-driven method to identify CVs in form of membership functions that span an invariant subspace of the Koopman operator associated with a stochastic dynamics system \cite{ISOKANN}.  

\subsection{Theoretical Background of ISOKANN}

In order to introduce the general idea behind ISOKANN, let us again consider the generator $\cL$ of the dynamical process under consideration and the associated Koopman operator $K_\tau =\exp(\tau \cL)$. Assume that the subspace $\cS_m=\textrm{span}\{\varphi_0,\ldots,\varphi_m\}$ of the dominant eigenfunctions of $K_\tau$ can be spanned by a partition of unity $\chi_0,\ldots,\chi_m$, i.e.,
\[
\cS_m=\mathrm{span}\{\chi_0,\ldots,\chi_m\}, \quad \sum_{i=0}^m \chi_i = 1,\quad \chi_i\ge 0.
\]
ISOKANN aims to compute such a partition of unity $\chi_0,\ldots,\chi_m$. 

Let $\chi(x)\in\R^{(m+1)}$ denote the vector $\chi(x)=(\chi_0(x),\ldots,\chi_{m}(x))^\top$, and let $\varphi(x)$ denote the same for the eigenfunctions $\varphi_i$.
Then $\cL\varphi=L\varphi$ with  $((\cL\varphi_i)(x)) = (L\cdot \varphi(x))_i$,
where $L$ denotes the $(m+1)\times (m+1)$ diagonal matrix that contains the eigenvalues $\lambda_0,\ldots,\lambda_m$ of $\cL$, $L_{ij}=\lambda_i\delta_{ij}$. Since both, the $\chi_i$ and the $\varphi_i$, span the subspace $\cS_m$, there is an invertible matrix $C\in\R^{(m+1)\times (m+1)}$ such that, for each $x\in\X$,
\[
\chi(x)=C\varphi(x),\quad\text{and}\quad \varphi(x)=C^{-1}\chi(x).
\]
Consequently, the matrix $\cQ= CL C^{-1}\in\R^{(m+1)\times (m+1)}$ satisfies
\begin{equation}\label{cLQ}
\cL\chi=\cQ\chi.
\end{equation}
We can compute the eigenvalues and eigenvectors of $\cQ$: Let $e_i$ denote the $i$-th unit vector containing $1$ at entry $i$ and otherwise $0$. Then we see that $v_i^\top=e_i^\top C^{-1}$ is the left eigenvector of $\cQ$ for eigenvalue $\lambda_i$,
\[
v_i^\top\cQ = \lambda_i v_i^\top,
\]
such that $\cQ$ inherits the dominant eigenvalues $\lambda_0,\ldots,\lambda_m$ of $\cL$.

\subsection{Effective Dynamics for ISOKANN}
By choosing $\chi$ as the CV, the effective dynamics in the ISOKANN latent space $\chi(\X)$ is given in terms of the $(m+1)$-dimensional column vector $z=\chi(x)$, and the projection $\Pi_\xi$ is averaging over the level sets of $\chi$. For ISOKANN, we can use (\ref{cLQ}) to see that the effective dynamics (\ref{eff-dynamics}) has the form
\begin{equation}\label{eff-dyn-chi}
dz_t = \cQ z_t\dt + \hsigma(z_t)\dW_t,
\end{equation}
with $z$-dependent $(m+1)\times (m+1)$ noise intensity matrix $\hsigma(z)=\sigma A(z)$ with $A$ from
\[
(AA^\top)_{lk}(z) = \Pi_\xi\Big(
    \sum_{i,j=1}^{d}  \frac{\partial \chi_l}{\partial x_i} \frac{\partial \chi_k}{\partial x_j}
    \Big)(z),
    \]
or from (\ref{effective-coeff-alternative}) with $\xi=\chi$. 
The resulting diffusion tensor is
\begin{equation}\label{eq:DiffusionTensor}
\begin{aligned}
\tilde{D}(z) 
&= \frac{1}{2} \sigma^2 A(z)A(z)^\top 
= \frac{1}{2} \hsigma(z)\hsigma(z)^\top.
\end{aligned}
\end{equation}
In general, the diffusion tensor is positive semidefinite since the dynamics on the latent space, i.e. the full simplex, are constrained by $\sum_0^m z_l = 1$, implying that one direction does not contribute to diffusion.

Following (\ref{eff-gen}), the associated effective generator has the form
\begin{equation}\label{eff-gen2}
\tL f(z)=(\cQ z)^\top\cdot \nabla_z f(z) + \frac{1}{2}\text{Tr}\left(\hsigma\hsigma^\top(z)\nabla_z^2 f(z)\right),
\end{equation}
where $\text{Tr}$ denote the trace operator.
Using the left eigenvectors $v_i$ of $\cQ$, we see that
\begin{equation}\label{ev:gen-eff}
\tphi_i(z)=v_i^\top z\quad \text{leads\ to}\quad \tL \tphi_i=\tilde{\lambda}_i\tphi_i,
\end{equation}
that is, $\tL$ inherits the dominant eigenvalues of $\cL$, and the associated eigenfunctions $\tphi_i(z)=v_i^\top z$ are linear in $z$. 
\paragraph{Dimension reduction $m+1\to m$:} $\chi$ contains  $m+1$ different functions $\chi_0,\ldots,\chi_m$. Since they form a partition of unity, we have that $\chi_m = 1 - \sum_{i=0}^{m-1} \chi_i$ and can thus reduce by one dimension. Let $\chi' = (\chi_0, \ldots, \chi_{m-1})$ be the reduced $m$-dimensional vector of the $m$ first $\chi$-functions, and let the matrix $\cQ$ be written in  block form:
\[
\cQ = \begin{bmatrix} \cQ' & b \\ a^\top & c \end{bmatrix}
\]
where $\cQ' \in \mathbb{R}^{m \times m}$, $a,b \in \mathbb{R}^{m \times 1}$, and $c \in \mathbb{R}$. Denoting with $\mathbf{1}\in\mathbb{R}^{m \times 1}$ the vector containing only ones, this yields
\begin{eqnarray*}
 \cQ  \chi = \begin{bmatrix}
(\cQ'-b\mathbf{1}^\top)\chi' + b\\
(a^\top-c\mathbf{1}^\top) \chi' + c
\end{bmatrix}
\end{eqnarray*}
such that the $m$-dimensional form of the effective dynamics (\ref{eff-dyn-chi}) reads
\begin{equation}\label{eff-dyn-chi-2}
dz'_t = b\dt + Q z'_t\dt + \hsigma(z'_t)\dW'_t,
\end{equation}
where $Q=\cQ'-b\mathbf{1}^\top$, $W'_t$ an $m$-dimensional vector of independent standard Wiener processes, and $\hsigma$ given by the same formula as above (just one dimension smaller). Thus, the reduced effective dynamics of ISOKANN is an Ornstein-Uhlenbeck process with linear drift and state-dependent noise intensity.

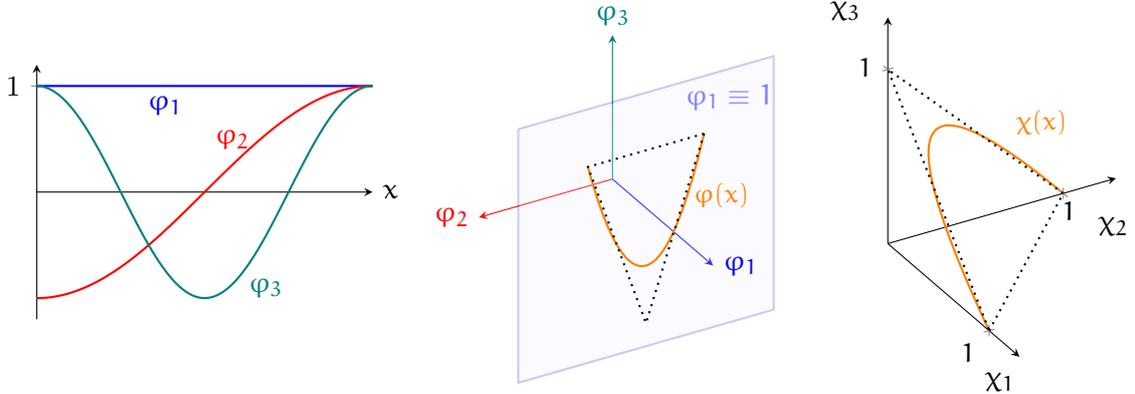
\begin{figure}[htbp]
\centering
% \documentclass{standalone}
% \usepackage{pgfplots}
% \pgfplotsset{compat=1.18}
% \usepgfplotslibrary{groupplots}

% \begin{document}
\begin{tikzpicture}
\begin{groupplot}[
  group style={
    group size=3 by 1,
    horizontal sep=2.3cm,
  },
  width=6cm, height=6cm,
  scale mode=scale uniformly,
  %scale only axis,
]

% ------------------------
% 1️⃣ Left: line plot φ_i(x)
% ------------------------
\nextgroupplot[
    axis lines=middle,
    axis line style={-stealth},
    xmin=0, xmax=pi,
    ymin=-1.2, ymax=1.2,
    xtick=\empty,
    ytick={0,1},
    yticklabels={0,1},
    xlabel={$x$},
    every axis x label/.style={at={(current axis.right of origin)}, anchor=west},
    every axis y label/.style={at={(current axis.above origin)}, anchor=south},
    clip=false,
  ]

  % φ1: linear
  \addplot[thick, blue, domain=0:pi, samples=100] {1};
  \node[anchor=west, blue] at (axis cs:pi/10*3,0.82) {$\varphi_1$};

  % φ2: quadratic
  \addplot[thick, red, domain=0:pi, samples=100] {sin(deg(x)-90)};
  \node[anchor=west, red] at (axis cs:pi/10*5,0.5) {$\varphi_2$};

  % φ3: sine-like
  \addplot[thick, teal, domain=0:pi, samples=100] {sin(deg(2*x)+90)};
  \node[anchor=west, teal] at (axis cs:pi/10*6, -.9) {$\varphi_3$};

% ------------------------
% 2️⃣ Middle: plane x = 1
% ------------------------
\nextgroupplot[
  view={60}{30},
  axis equal image,
  axis lines=middle,
  xshift=-9mm,
  yshift=3mm,
  width=0.9*6cm,               % 90 % of the original 6 cm
  height=0.9*6cm,              % 90 % of the original 6 cm
  scale only axis=false,       % important! scales the entire plot including arrows/labels
  % … all your other options (view, axis lines=left, etc.) …
  %axis line style={->},
  xlabel={$\varphi_1$}, xlabel style={text=blue},
ylabel={$\varphi_2$}, ylabel style={text=red},
zlabel={$\varphi_3$}, zlabel style={text=teal},
  xmin=0, xmax=3,
  ymin=0, ymax=2.3,
  zmin=0, zmax=2.5,
  xtick={0,1},
  ytick={0,1},
  ztick={0,1},
  ticks=none,
  xticklabels={},
  yticklabels={},
  zticklabels={},
  x axis line style={blue},    % or any color you want
  y axis line style={red},
  z axis line style={teal},
  %tick style={draw=black, major tick length=3pt},
  %grid=none,
  samples y=1, % dont close the curve
  y dir=reverse,
  xlabel style={at={(xticklabel cs:1)}, anchor=west},   
  ylabel style={at={(yticklabel cs:-0)}, anchor=east}, 
  zlabel style={at={(zticklabel cs:1)}, anchor=south},
  clip=false
]

% Plane
\addplot3[
  fill=blue!10,
  draw=blue!70!black,
  opacity=0.2,
  thick,
]
coordinates {
  (1,2.2,2)
  (1,-2.2,2)
  (1,-2.2,-2.4)
  (1,2.2,-2.4)
  (1,2.2,2)
};
\node[above right, text=blue!50] 
    at (axis cs:1, -.5, 1.4) {$\varphi_1 \equiv 1$};

\addplot3[
  thick,
  color=orange,
  samples=100,
  domain=0:pi,
]
({1}, {sin(deg(x)-90)}, {sin(deg(2*x)+90)});

\addplot3[dotted, thick]
coordinates {
    (1, -1, 1)
    (1, 0, -2)
    (1, 1, 1)
     (1, -1, 1)
};

\node[orange, above, font=\small\bfseries, inner sep=2pt] 
    at (axis cs:1, -1.3, -.5) {$\varphi(x)$};

% fixed point of koopman dynamics
% \addplot3 [
%     scatter,           % activates scatter mode
%     only marks,        % plots nothing except the marks
%     mark=x ,           % filled circle (you can also use mark=ball, mark=square*, etc.)
%     color=black,
%     mark size=2pt,     % adjust size as you like
%     % or black, red, whatever you prefer
% ] coordinates {(1,0,0)};

% ------------------------
% 3️⃣ Right: simplex
% ------------------------
\nextgroupplot[
  yshift=2mm,
  view={60}{30},
  axis equal image,
  axis lines=middle,
  xlabel={$\quad\chi_1$}, ylabel={$\chi_2$}, zlabel={$\chi_3$},
  xmin=0, xmax=1.3,
  ymin=0, ymax=1.3,
  zmin=0, zmax=1.3,
  xtick={0,1}, ytick={0,1}, ztick={0,1},
  ticklabel style={font=\small},
  grid=none,
  ticks=both,
  minor tick num=0,
  declare function={
    sA(\x) = 1;
    sB(\x) = sin(deg(\x) - 90);    % sin(x - 90°)
    sC(\x) = sin(deg(2*\x) + 90);  % sin(2x + 90°)
  },
  xlabel style={at={(xticklabel cs:1)}, anchor=south},
    ylabel style={at={(yticklabel cs:1)}, anchor=east},
    zlabel style={at={(zticklabel cs:0.9)}, anchor=south},
    scale only axis=true,
]
% Scatter
\addplot3[
    thick,
    color=orange,
    samples=50,
    domain=0:pi,
    samples y=1,
]
(
    {(0.33 * sA(x)  - 0.5 * sB(x) + 0.17 * sC(x))},
    {(0.33 * sA(x)  + 0.5 * sB(x) + 0.17 * sC(x))},
    {(0.33 * sA(x)  + 0   * sB(x) - 0.33 * sC(x))}
);

\node[orange, above, font=\small\bfseries, inner sep=2pt] 
    at (axis cs:1, .3, 1) {$\chi(x)$};

% Simplex edges
% \addplot3[thick, color=black] coordinates {(1,0,0) (0,1,0)};
% \addplot3[thick, color=black] coordinates {(0,1,0) (0,0,1)};
% \addplot3[thick, color=black] coordinates {(0,0,1) (1,0,0)};

\addplot3[dotted, thick]
coordinates {
    (1, 0, 0)
    (0, 1, 0)
    (0, 0, 1)
    (1, 0, 0)
};

\end{groupplot}
\end{tikzpicture}

% \end{document}
\caption{Illustration of the transformation from eigen- to $\chi$-functions.
Left: Dominant eigenfunctions $\varphi_i(x)$. Middle: 
Eigenfunction components plotted in the plane $\varphi_1 \equiv 1$, with an enclosing simplex indicated.
Right: The transformation that maps this simplex to the unit simplex sends $\varphi$ to $\chi$, which then lies on the unit simplex.}
\label{fig:koopman}
\end{figure}

%\AS{I tried to emphasize the geometric side of this transform}
This reduction of dimension from $m+1 \rightarrow m$, or equivalently from working with $z$ to $z'$, reflects the fact that the (m+1)-tuple $z$ lies on the unit m-simplex embedded in $\mathbb{R}^{m+1}$,
%\AS{@Luca, would it be a problem changing the indexing to l=1..m+1 ?}
\begin{equation}
\bbDelta^{m} = \left\lbrace z\in[0,1]^{m+1} : \sum_{l=0}^m z_l = 1 \right\rbrace,
\label{full-simplex}
\end{equation}
and therefore contains one redundant coordinate.
Eliminating this redundancy amounts to projecting $\bbDelta^m$ onto an intrinsic $m$-dimensional representation in $\mathbb{R}^m$, yielding the reduced simplex
\begin{equation}
    \bar{\bbDelta}^{m} = \left\lbrace z'\in[0,1]^m : \sum_{l=0}^{m-1} z'_l \le 1 \right\rbrace.
    \label{reduced-simplex}
\end{equation}
Equation (\ref{eff-dyn-chi-2}) therefore represents the projection of (\ref{eff-dyn-chi}) onto this lower-dimensional simplex, see Fig.~\ref{fig:koopman}.

One advantage of this perspective is that the diffusion matrix becomes strictly positive definite: $\tilde{D}(z') \succ 0$.
%
%\AS{I think this last paragraph is not necessary:}
%
Indeed the Jacobian matrix $\nabla \chi$ of the CV $\chi$ has full rank $m$ on the reduced simplex and the noise $\sigma$ acts in all directions in full dynamics.

%%%%%%%%%%%%%%%%%%%%%%%%%%%%%%%%%%%%%%%%
%%%%%%%%%%%%%%%%%%%%%%%%%%%%%%%%%%%%%%%%
\subsection{The ISOKANN algorithm}

\newcommand{\KK}{K_\tau}

Without sufficient sampling of rare but essential transitions between metastable states, any method approximating eigenfunctions will struggle. Therefore, in \cite{sikorski2024learning} optimal importance sampling was integrated into ISOKANN, yielding a dramatic variance reduction when estimating the action of the Koopman operator. In the following, we focus on the algorithmic application of ISOKANN in higher-dimensional cases by incorporating the Inner Simplex Algorithm, which enables the computation of multi-dimensional $\chi$ functions.

In order to compute the invariant subspace $\cS_m$ in terms of the vector valued $\chi$ function, ISOKANN combines three different techniques: an iterative procedure similar to the power iteration, representation of $\chi$ through neural networks and data-driven Monte-Carlo approximation of the action of $\KK$.

Let us first look at the iteration that lies at the core of the method.
Similarly to the power iteration, the Arnoldi iteration, or more general Krylov-subspace methods, the underlying idea is that the iterative application of $\KK$ to any generic initial function will converge to the dominant eigenfunction, 
\[ \KK^n b_0 \rightarrow c \varphi_0,\quad b_0 \in \cH, c\in \mathbb{R},\]
which is a direct consequence of the eigenexpansion (\ref{eq:eigenexpansion}).

Since we are not interested in the dominant eigenfunction $\varphi_0\equiv 1$ itself, but in the following slowest decaying eigenfunctions, we modify this power iteration scheme by constructing a sequence of linear operators $A_n$ such that the iteration $\chi^{(n+1)}=A_n\chi^{(n)}$ converges to the dominant subspace of $\KK$.

To retain the exponentially decaying contributions from these slower eigenfunctions, we 
compose the application of the Koopman operator with a dynamically chosen linear transformation $S_n$,
\begin{align}
\label{eq:xSKx}
    \chi^{(n+1)} = S_n \KK \chi^{(n)},
\end{align}
where $\chi^{(n)}: \X \rightarrow \mathbb{R}^{m+1}$ denotes the current iterate and $S_n \in \mathbb{R}^{(m+1) \times (m+1)}$ is a matrix computed for each iteration such that it rescales the dominant modes to unit magnitude, ensuring stability and preventing collapse to the trivial eigenfunction.
This construction allows to retain the $m+1$ dimensional dominant subspace whereas all faster modes keep decaying, ensuring convergence 
of $\cS_m^{(n)} = \text{span} \{\chi^{(n)}_1,\ldots,\chi_{m+1}^{(n)}\}$ 
to $\cS_m$.

\newcommand{\bK}{\mathbf{K}}
Assuming convergence, i.e. $\chi^{(n)} \approx S_n \KK \chi^{(n)}$ for some large $n$, we see that the matrix $\bK^{(n)}=S_n\inv$ is the finite dimensional representation of $\KK$ on the dominant subspace $\mathcal{S}_m^{(n)}$, $\bK^{(n)} \chi \approx \KK \chi$. 
Since $\KK =\exp(\tau \cL)$ we also have $\bK^{(n)} \approx \exp(\tau \cQ )$ which in turn allows us to obtain the rate matrix $\cQ$ directly from $\bK^{(n)}$ through the matrix logarithm.

Since $\mathcal{X}$ can be high-dimensional, we approximate $\chi^{(n)}$ using a neural network $\chi_{\theta_n}$ with trainable parameters $\theta_n$. Neural networks can represent smooth low-dimensional manifolds embedded in high-dimensional space, which is particularly suited for this task of learning reaction coordinates with low intrinsic dimension.
However, it also means that the equality in (\ref{eq:xSKx}) must be realized through a training task. Here we choose the classical regression by minimizing the residual in terms of the mean squared error (MSE) over a set of data points $x_i\in \X,\;i = 1,\ldots, N$. 
The loss function then reads
\begin{align}\label{eq:loss}
    L(\theta_{n+1}) = \frac{1}{N}\sum_{i=1}^N \left\| \chi_{\theta_{n+1}}(x_i) - T_n(x_i)\right\|^2,\quad\text{with}\quad T_n(x_i) = S_n \KK \chi_{\theta_n}(x_i)
\end{align}
where the target $T_n$ is based on the previous iteration and as such is independent of the current iteration's weights (i.e. with "frozen gradients" w.r.t. the parameter optimization for $\theta_{n+1}$).

The evaluation of the Koopman operator $\KK \chi^{(n)}(x_i)=\KK\chi_{\theta_n}(x_i)$ is realized through the Monte-Carlo approximation of its expectation value formulation (\ref{eq:defkoopman}):
\[
\KK \chi^{(n)}(x_i) = \frac{1}{M}\sum_{k=1}^M \chi^{(n)}(y_{i,k})\]
where $y_{i,k} \sim P(X_\tau | X_0 = x_i)$ are $M$ independent samples that evolved in time $\tau$ from their starting positions $x_i$.
In this regard, the training data for ISOKANN consists of two tensors: 
\[
X \in \mathbb{R}^{d \times N}, \quad 
Y \in \mathbb{R}^{d \times M \times N}
\]
where the columns $x_i = X_{:,i} \in \mathbb{R}^d$ are the $N$ initial points of the short burst simulations and $y_{i,k}=Y_{:,k,i} \in \mathbb{R}^d$ the endpoints of the respective $M$ bursts.

\paragraph{Estimation of $S_n$ through the inner simplex algorithm:}
Let us for a moment consider the ISOKANN iteration (\ref{eq:xSKx}) without the regularizing effect of $S_n$, i.e. with $S_n = \text{Id}$.
In particular, consider the image of the data points: 
$\KK\chi^{(n)} (X) 
=\left\{\KK \chi^{(n)}(x_i) \mid x_i \in X \right\} 
%= \left\{\KK^n \chi_0(x_i) \mid x_i \in X \right\} 
\subset \mathbb{R}^{m+1}$.
Due to the spectral properties of $\KK$ ($\Lambda \le 1$, single eigenvalue 1) this set will, for $n\rightarrow\infty$, converge to a unique fixed point in $\mathbb{R}^{m+1}$.

The \emph{inner simplex algorithm} (ISA) \cite{WeberGalliat2002} counteracts this contraction by identifying the most extremal points of each iterate $\KK \chi^{(n)}(X)$ and constructing the linear map $S_n$ that sends these points to the vertices of the unit simplex in $\mathbb{R}^{m+1}$.
In each iteration, the extremal points are selected recursively by maximal distance to the span of all previously chosen points, starting with the point of largest magnitude.
Thus, ISA determines a simplex inscribed in the convex hull of the data, successively expanding $(m+1)$ independent directions to unit scale and thereby stabilizing the iteration (\ref{eq:xSKx}).
Because the slowest modes will dominate the extrema, their contributions prevail under iteration, while faster modes decay.
Note that unlike the Koopman operator, which is evaluated at each data point $x_i$ individually, the estimation of $S_n$ in ISA is computed over the entire set $\KK \chi^{(n)}(X)$; once determined, the resulting matrix $S_n$ can then be applied to each vector in $\KK \chi^{(n)}(X)$. 

The ISA generalizes the shift–scale operation used in earlier ISOKANN works \cite{ISOKANN,sikorski2024learning,Donati2024,SikorskiRabbenChewleWeber,Donati2025} from the one-dimensional case -- where $\chi$ represented a 1-simplex in $\mathbb{R}^2$ and $S_n$ was an explicit affine-linear map -- to the higher-dimensional setting.
Conceptually, ISA can be viewed as a higher-dimensional analogue of this shift–scale step, and itself is a subroutine of PCCA+ \cite{Deuflhard2004, Kube2007, Weber2018}, the method that originally inspired the idea of $\chi$-functions.

\begin{algorithm}[H]
\caption{ISOKANN Algorithm for $m+1$ dimensions}
\label{alg:isokann}

\SetAlgoLined
\KwIn{
    Initial features $X \in \mathbb{R}^{d \times N}$ with columns denoted $x_i$;\\
    Propagated samples $Y \in \mathbb{R}^{d \times M \times N}$ with columns denoted $y_{i,k}$;\\
    Neural network $\chi_\theta:\mathbb{R}^d \to \mathbb{R}^{m+1}$
}
\KwOut{
    Trained network $\chi_\theta$ approximating $\mathcal{S}_m$
}

\BlankLine
Initialize parameters $\theta_0$ randomly. Set $n \gets 0$.

\Repeat{convergence}{

  \tcc{Monte Carlo Koopman evaluation}
  $K_i \gets \frac{1}{M} \sum_{k=1}^M \chi_{\theta_n}(y_{i,k}) \quad \forall i=1,\dots,N$\;

  \tcc{Inner Simplex Algorithm to construct $S$}
  \Begin{  % start indent
    $j_1 \gets \arg \max_k \| K_k\|$\;
    \For{$i = 1$ \KwTo $m$}{
          $j_{i+1} \gets \arg\max_k \operatorname{dist}\!\left(K_k, 
          \operatorname{span}\{K_{j_1},\dots,K_{j_{i}}\}
          \right)$\;
    }
    $S \gets [\,K_{j_1}\;\; K_{j_2}\;\; \cdots\;\; K_{j_{m+1}}\,]^{-1}$\;
  }  % end indent

  \tcc{Form training targets}
  $T_i \gets S\,K_i \quad \forall i=1,\dots,N$\;

  \tcc{Mean-squared training loss}
  $L(\theta) \gets \frac{1}{N} \sum_{i=1}^N \|\chi_\theta(x_i) - T_i\|^2$\;

  \tcc{Update network parameters}
  $\theta_{n+1} \gets \arg\min_\theta L(\theta)$\;

  $n \gets n+1$\;
}

\Return $\chi_\theta$.

\end{algorithm}

\paragraph{Algorithmic considerations} We now discuss several practical aspects relevant to implementation.

\begin{itemize}
\item \textbf{Initialization:}
For the power-iteration scheme to work, the initial iterate $\chi_0$ must contain nonzero contributions from the dominant modes. In practice, random initialization of the neural network parameters satisfies this condition almost surely.
    
\item \textbf{Architecture:}
In ISOKANN, the neural network $\chi_\theta:\mathbb{R}^d \to \mathbb{R}^{m+1}$ acts as a black-box regressor; thus, the specific architecture is not crucial. We typically employ fully connected feed-forward networks with 2–4 hidden layers and sigmoid or ReLU activations. For sigmoid activations, the final layer uses the identity function, allowing the output to map onto $[0,1]$.

\item \textbf{Features:} 
A suitable choice of features can greatly simplify learning. Accordingly, $X$ and $Y$ may represent physically meaningful features (e.g., normalized pairwise particle distances) rather than raw coordinates.

\item \textbf{Network reuse:}
Although in theory each ISOKANN iterate $\chi^{(n)}$ is represented by a separate network $\chi_{\theta_n}$, in practice we continue training a single model, updating its weights iteratively. This approach treats each iteration as a refinement rather than a full retraining.

\item \textbf{Training loop:}
The training consists of an inner loop of stochastic gradient descent over multiple epochs and potentially minibatches, nested inside the outer ISOKANN iteration. Reusing network weights allows the number of training epochs per ISOKANN iteration to be reduced to one, ensuring that each epoch trains on the most recent target update.

\item \textbf{Inner simplex algorithm:}
The ISA constructs the matrix $S_n$ mapping the extremal points of $\KK \chi^{(n)}$ onto a unit simplex and thereby defines the training targets $T$. Since the orientation of the simplex (i.e., vertex permutation) is arbitrary, it is necessary to permute the components of $T$ to align with the previous $\chi^{(n)}$, for example by maximizing the diagonal of their covariance.

\item \textbf{Optimizer:}
We employ either Adam or Nesterov Momentum (learning rate $10^{-3}$), observing that the latter can outperform Adam when the optimization stagnates. To prevent overfitting -- which is especially important because the network is evaluated at unseen positions during the Koopman estimation -- we add Tikhonov regularization with decay rate $10^{-4}$.

\item \textbf{Convergence:}
Assessing convergence is nontrivial; in practice, we use a fixed number of outer iterations informed by prior experience, or alternatively monitor empirical or validation loss values.

\item \textbf{Validation loss:}
The validation loss is bounded below by the Monte-Carlo noise level in the Koopman estimation. Thus, it may be advantageous to trade a large number of validation points for fewer, higher-fidelity samples (i.e., larger number $M$ of $y_{i,k}$ per $x_i$). In this case, $S$ should still be computed on the training data, since sparse validation points may not span the full range of $\chi$. This setup can lead to the seemingly paradoxical situation where the empirical loss stagnates when approaching the data noise level while the validation loss keeps decreasing indicating further improvement.

\item \textbf{Initial samples $X$:}
The initial data points $x_i$ used for training do not need to be drawn from the stationary distribution of the dynamics. Their distribution mainly affects the weighting of residuals in the loss function rather than the target mapping itself. Consequently, ISOKANN remains robust when trained on non-stationary data and readily supports adaptive sampling strategies, e.g. where new points are generated based on the evolving reaction coordinate $\chi^{(n)}$.

\item \textbf{Koopman samples $Y$:}
For a single training point $x_i$ with $M$ Koopman samples, the losses
\begin{eqnarray}
L_1(\theta) & = & \left\| \chi_\theta(x_i) -\frac{1}{M}\sum_{k=1}^M S\chi_\tn (y_{i,k})\right\|^2 \nonumber\\
& = & \chi_\theta(x_i)^2 - \frac{2}{M}\chi_\theta(x_i)\sum_{k=1}^M S\chi_\tn(y_{i,k})+\frac{1}{M^2}\left(\sum_{k=1}^M S\chi_\tn(y_{i,k})\right)^2 \\
L_2(\theta) & = &\frac{1}{M}\sum_{k=1}^M \left\| \chi_\theta(x_i) -S\chi_\tn(y_{i,k})\right\|^2 \nonumber\\
& = & 
    \chi_\theta(x_i)^2 - 
    \frac{2}{M}\chi_\theta(x_i)\sum_{k=1}^MS\chi_\tn(y_{i,k}) + \frac{1}{M}\sum_{k=1}^M \left(S\chi_\tn(y_{i,k})\right)^2
\end{eqnarray}
yield identical gradients, $\nabla_\theta L_1 = \nabla_\theta L_2$. Thus, averaging over multiple Koopman samples can equivalently be realized by multiple replicas of single Koopman samples.

Increasing $M$ reduces variance only as $O(1/M)$ while increasing computational cost linearly. In practice, much of the averaging effect may already be achieved implicitly through regularization by the network’s smoothness and the distribution of nearby training points: values of $\chi_\theta$ at close states tend to co-vary smoothly, approximating a local ensemble average.

Consequently, rather than increasing $M$, it can be more effective to sample a larger number of distinct initial states $x_i$ with fewer Koopman samples each -- down to $M=1$ in the extreme case -- trading fidelity for broader coverage of the state space.

\item \textbf{Single trajectory data:}
Since $X$ does not need to sample the stationary distribution, it is possible to restrict attention to a subset of the state space, effectively analyzing the dynamics on this subset. Using the $M=1$ argument above, one might consider using a single trajectory as $x_i$ and its time-shifted points as $y_i$. However, this data only reflects the transport of $\chi$ values along the trajectory. In this setting, we recommend sufficient regularization and, for reversible processes, augmenting the trajectory with its time-reversed counterpart. With these precautions, ISOKANN can extract the slowest process from a single simulation trajectory \cite{SikorskiRabbenChewleWeber}.

\end{itemize}

%%%%%%%%%%%%%%%%%%%%%%%%%%%%%%%%%%%%%%%%%%%%%%%%%%%%%%%%%%%%%%%%%
%%%%%%%%%%%%%%%%%%%%%%%%%%%%%%%%%%%%%%%%%%%%%%%%%%%%%%%%%%%%%
\section{Transition Rates and Pathways}\label{sec:transitions}

There is extensive literature on computing the transition rates and pathways in diffusive or molecular processes, ranging from TPS \cite{TPS}, interface methods like milestoning \cite{Elber:04} or TIS \cite{TIS}, or \emph{transition path theory} (TPT) \cite{TPT0,TPT1,TPT2}, see \cite{Szabo2019} for a conceptual review. We will herein concentrate on TPT and how it can be used to investigate mean first passage times, committor functions, and dominant transition pathways in full and effective dynamics.  

\subsection{Transition kinetics of the full dynamics}

While the slow modes of the full dynamics are given by the dominant eigenfunctions of its Koopman operator $K_\tau$ and its generator $\cL$, we can also understand the transition kinetics into a selected set $A\subseteq\X$ with a smooth boundary or between $A$ and a disjoint set $B\subseteq \X$. In order to provide a short review, we first introduce the first hitting time $\tau_x(B)$ of set $B$ if started in $x\in\X$, which is a random time  formally defined by $\tau_x(B)=\inf_{t\ge 0}\{X_t\in B|X_0=x\}$. The \emph{mean first passage time} (MFPT) from $x$ to $B$, $m_B(x)=\mathbb{E}(\tau_x(B))$, is given by the boundary value problem
\begin{equation}\label{MFPT}
\cL m_B(x)=-1,\;x\in\X\setminus B\quad \text{and}\quad m_B(x)=0,\;x\in B.
\end{equation}
The \emph{committor function} $q(x)$ is the probability to reach $B$ before $A$ when starting from $x$:
\begin{equation}
q(x) = \mathbb{P}(\tau_x(B) < \tau_x(A)),\quad\text{for}\quad x\in \X\setminus(A\cup B),
\end{equation}
and can also be characterized as the unique solution of a boundary value problem,
\begin{equation}\label{committor}
\cL q(x) = 0 \;\text{for}\;  x\in \X\setminus(A\cup B),\quad\text{and}\quad 
q(x)=0, \;\text{for}\; x\in A,\quad
q(x)=1, \;\text{for}\; x\in B.
\end{equation}
TPT describes the ensemble of reactive trajectories \cite{TPT2,TPTillu}— those paths that go from 
$A$ to $B$ without returning. The \emph{reactive density}
\[
\mu_{AB}(x)=\mu(x)q(x)(1-q(x)),
\]
is the stationary probability density of finding the system while on a reactive trajectory (for reversible systems).
For processes of the form (\ref{eq:full}), the \emph{reactive probability flux} associated with transitions from $A$ to $B$ is
\begin{equation}\label{prob_flux}
j_{AB}(x)=\frac{1}{2}\mu(x)\,\sigma^2\,\nabla q(x),
\end{equation}
where $\mu(x)\propto \exp(-\beta V(x))$ again denotes the invariant density of the process, and $\mu_A=\int_A \mu(x) \, dx$ denotes the stationary density of the set $A$.
The \emph{reactive transition rate} from $A$ to $B$ is given by the total flux through any surface $S$ separating $A$ and $B$,
\begin{equation}
k_{AB} = \frac{1}{\mu_A}\int_S j_{AB}(x)\cdot n_S(x)\,dS(x),
\end{equation}
where $n_S$ is the unit normal to any dividing surface $S$ separating $A$ from $B$. In the reversible case, this flux is independent of the chosen dividing surface.
Equivalently, one can express the rate via the Dirichlet form:
\begin{equation}\label{eq:full_rate}
k_{AB} = \frac{1}{\mu_A} \frac{\sigma^2}{2}\int_{\X\setminus(A\cup B)} 
\mu(x)\,\nabla q(x)^\top\,\nabla q(x)\,dx.
\end{equation}

The \emph{transition pathways} are extracted from the reactive flux field $j_{AB}$: 
Its streamlines (integral curves) give the most probable reaction pathways.
These can be visualized as “reaction tubes” showing how probability flows between $A$ and $B$;
high-magnitude regions of $j_{AB}$ identify the dominant reaction tubes or channels and bottlenecks of the transition.

\subsection{Transition kinetics of the effective dynamics}
In what follows, we consider the computation of kinetic quantities for the effective dynamics after dimension reduction, such as MFPT, committor, reactive flux and transition rate for disjoint sets $\tilde{A},\tilde{B}$ of the latent space $\mathbb{Z}$ which is identified with the reduced simplex $\bar{\bbDelta}^{m}$ (\ref{reduced-simplex}).
In (\ref{eff-dyn-chi-2}) the reduced coordinates on the simplex 
were denoted by $z'$ to emphasize the projection from the full 
simplex.
From this point onward, for simplicity of notation, we omit the prime and continue to denote all quantities (variables, functions and operators) using the same symbols as for the effective dynamics. 
%
% We can compute the kinetic quantities, MFPT, committor, reactive flux and transition rate, for the effective dynamics for disjoint sets $\tilde{A},\tilde{B}$ of the latent space $\mathbb{Z}$. However, in order to compare these quantities to those of the full dynamics, we have to make preparations: 

When the sets $\tilde{A},\tilde{B}\subseteq \bar{\bbDelta}^{m}$ for the effective dynamics are chosen, we have to define the sets $A,B$ with respect to them,
\[
A=\chi^{-1}(\tilde{A}),\quad\text{and}\quad B=\chi^{-1}(\tilde{B}),
\]
assuming that $\chi$ is invertible. While the formula for the effective MFPT and the effective committor $\tilde{q}$ are still given by (\ref{MFPT}) and (\ref{committor}) by replacing $\cL$, $x$, $A$, and $B$ by $\tL$, $z$, $\tilde{A}$, and $\tilde{B}$, the equation for the reactive flux and the transition rate now read
\begin{eqnarray}
\tilde{j}_{AB}(z) & = & \frac{1}{2} \tilde{\mu}(z)\,\hsigma\hsigma^\top(z)\,\nabla \tilde{q}(z)\label{eff-flux}\\
\tilde{k}_{AB} & = & \frac{1}{\tilde{\mu}_{\tilde{A}}}
\frac{1}{2}\int_{\mathbb{Z}\setminus(\tilde{A}\cup \tilde{B})} 
\tilde{\mu}(z)\,\nabla \tilde{q}(z)^{\top} \hsigma\hsigma^\top(z)\nabla \tilde{q}(z)\,dz\label{eff-rate},
\end{eqnarray}
where $\tilde{\mu}$ denotes the invariant density of the effective dynamics, and $\tilde{\mu}_{\tilde{A}}=\int_{\tilde{A}} \tilde{\mu}(z)\, dz$.
%%%%
%

In the case $m=1$, after reduction of the effective dynamics to dimension $1$, the generator is scalar in latent space $\mathbb{Z}=[0,1]$,
\[
\tL f(z)=(a+\lambda z) \nabla_z f(z) +  \frac{1}{2}\hsigma^2(z)\nabla_z^2 f(z),
\]
By defining the rate function $\phi$ and the \emph{effective potential} (or Pontryagin's potential going back to \cite{AndronovPontryaginVitt1933}), 
\begin{eqnarray}
\phi(z) & = & \int_{z^\ast}^z \frac{2(a + \lambda y)}{\hsigma(y)^2} \, dy,\label{phi}\\
V_{\chi}(z) & = & \log\Big(\frac{1}{2}\hsigma(z)^2\Big) - \phi(z)\,+\,\text{const},\label{eff_pot_1} 
\end{eqnarray}
the effective invariant density can be written
\[
\tilde{\mu}(z)=\frac{1}{Z_\chi} \exp(-V_\chi(z)),\qquad Z_\chi=\int_0^1 \exp(- V_\chi(z))dz.
\]
The potential $V_\chi$ still depends on the reference point $z^\ast$, 
but different choices of $z^\ast$ only add a constant to $V_\chi$, 
which cancels out in normalized quantities such as $\tilde{\mu}$.
%
%However,  different choices of $z^\ast$ only lead to  different additional $z$-independent constants in $V_\chi$ which cancel out in all final results. (2) $V_\chi$ depends nonlinearly on the inverse temperature $\beta=2/\sigma^2$ of the original dynamics. Therefore, the dependence of the effective invariant density on $\beta$ does not exhibit the standard form with $\beta$ as a linear scaling factor in the exponential.  
%

The MFPT of the effective dynamics from $z$ to set $\tB=[z_B,1]$ can be written down explicitly, 
\begin{eqnarray}
\tilde{m}_{\tB}(z) & = & \int_z^{z_B}\frac{2}{\hsigma(z)^2}e^{ V_\chi(z)}\,\int_{0}^{z} e^{- V_{\chi}(u)} \, du \,dz.\label{eff-MFPT-1d}
\end{eqnarray}
Furthermore, the effective committor to go from  set $\tA=[0,z_A]$ to set $\tB=[z_B,1]$ with $z_A<z_B$
can be written in the following form for $z\in [z_A,z_B]$:
\begin{equation}\label{eff-committor-1d}
\tilde{q}(z)=\frac{1}{Z_q}\int_{z_A}^z \frac{2}{\hsigma(z)^2} e^{V_\chi(z)} dz,\quad \mathrm{with}\quad Z_q =\int_{z_A}^{z_B} \frac{2}{\hsigma(z)^2} e^{ V_\chi(z)}dz.
\end{equation}
For the reactive flux and the reactive transition rate from $\tA$ to $\tB$ this leads to 
\begin{equation}\label{transition_rate+flux}
\tilde{j}_{AB}(z)=\tilde{k}_{AB} =\frac{1}{Z_q Z_\chi}.
\end{equation}
%
%%%%%%%%%%%%%%%%%%%%% m > 1 C A S E
%

In the case $m>1$, the effective generator (\ref{eff-gen}) on the reduced simplex $\bar{\bbDelta}^{m}$ reads
\begin{equation}
\label{eff_gen}
\begin{aligned}
\tL f(z) 
&= 
(b + Qz) \cdot \nabla_z f(z)
+ \frac{1}{2}
\text{Tr}
\left(
(\hsigma\hsigma^\top)(z)
\nabla_z^2 f(z)
\right),
\end{aligned}
\end{equation}
%
% \begin{equation}
% \label{eff_gen}
% \begin{aligned}
% \tL f(z) 
% &= 
% (b + Qz) \cdot \nabla_z f(z)
% + \frac{1}{2}
% \sum_{l,k=0}^{m-1}
% (\hsigma\hsigma^\top)_{kl}(z)
% \frac{\partial^2 f}{\partial z_l \partial z_k}(z) \, ,
% \end{aligned}
% \end{equation}
where $b$ and $Q$ were defined in (\ref{eff-dyn-chi-2}).
The corresponding effective potential is given by
\begin{equation}
\begin{aligned}
\label{eff_pot_m_dim}
V_{\chi}(z) 
&= 
\int_{\gamma:z^\ast\to z} \psi(y) \cdot dy,
\end{aligned}
\end{equation}
where
\begin{equation}
\begin{aligned}
\psi(z) 
&= 
2 (\hsigma\hsigma^\top)^{-1} 
\left(
\frac{1}{2}
\nabla_{z} \cdot (\hsigma\hsigma^\top)
-
(b+Qz)(z)
\right)
\end{aligned}
\end{equation}
and the line integral is taken along any smooth path $\gamma\subset \bar{\bbDelta}^{m}$ from a starting point $z_\ast$ to $z$. 
Since, under the reduction procedure, the diffusion matrix $\hsigma\hsigma^{\top}$ defined in (\ref{eq:DiffusionTensor}) is positive definite, it is invertible everywhere on $\bar{\bbDelta}^m$.
Moreover, the reversibility of the effective generator ensures that the integrand $\psi(z)$ is curl-free and that the effective potential is path-independent.
Note that if $m=1$, (\ref{eff_pot_m_dim}) reduces to (\ref{eff_pot_1}).

Finally, the effective invariant density is given by
\[
\tilde{\mu}(z)=\frac{1}{Z_\chi} \exp(- V_\chi(z)),\qquad Z_\chi=\int_{\bar{\bbDelta}^{m}} \exp(- V_\chi(z)) \, dz,
\]

Unlike the $m=1$ case, explicit expressions for the effective mean first-passage time $\tilde{m}_{\tilde{B}}$ and the effective committor $\tilde{q}$ cannot be derived. 
However, both quantities can be computed numerically by solving the corresponding boundary-value problems on $\bar{\bbDelta}^{m}$. 
Specifically, the effective MFPT satisfies
\[
\tL \tilde{m}_{\tilde{B}}(z) = -1, \quad z\in \bar{\bbDelta}^{m}\setminus\tilde{B},
\qquad
\tilde{m}_{\tilde{B}}(z)=0, \quad z\in\tilde{B},
\]
while the effective committor solves
\begin{equation}
\label{committor_m}
\tL \tilde{q}(z) = 0, \quad z\in \bar{\bbDelta}^{m}\setminus(\tilde{A}\cup\tilde{B}),
\qquad
\tilde{q}(z)=0 \text{ on } \tilde{A}, \quad
\tilde{q}(z)=1 \text{ on } \tilde{B}.
\end{equation}
%%%%%%%%%%%%%%%%%%%%%%%%%%%%%%%%%%%%%%%%%%%%
\section{Numerical Experiments}
\label{sec:numerical_experiments}
%
%In this section, we demonstrate the ISOKANN approach for computing collective variables (CVs) and effective dynamics on three benchmark systems of increasing complexity. 
%
%The objective is to validate the theoretical framework presented above by direct comparison between the full stochastic dynamics, its ISOKANN-based reduced representation, and the corresponding transition kinetics derived from transition path theory.
%
%In all three benchmark cases, we consider diffusive dynamics in full space as defined in (\ref{eq:full}).

In this section, we apply the ISOKANN approach to compute collective variables (CVs) and effective dynamics for three benchmark systems, each increasing in complexity. The primary goal is to validate the theoretical framework introduced in the previous sections by comparing the full stochastic dynamics with the ISOKANN-based reduced representation, and the corresponding transition kinetics derived from transition path theory.

For all three systems, we model diffusive dynamics in the full space given by overdamped Langevin dynamics, see (\ref{eq:full}). The ISOKANN algorithm is capable of constructing $\chi$-functions even in very high-dimensional state spaces. However, to facilitate direct comparison with ground truth, we limit our examples to systems of up to three dimensions, which allows us to discretize the full-dimensional dynamics and compare the reduced models directly.

\begin{figure}
    \centering
    \includegraphics[width=1\linewidth]{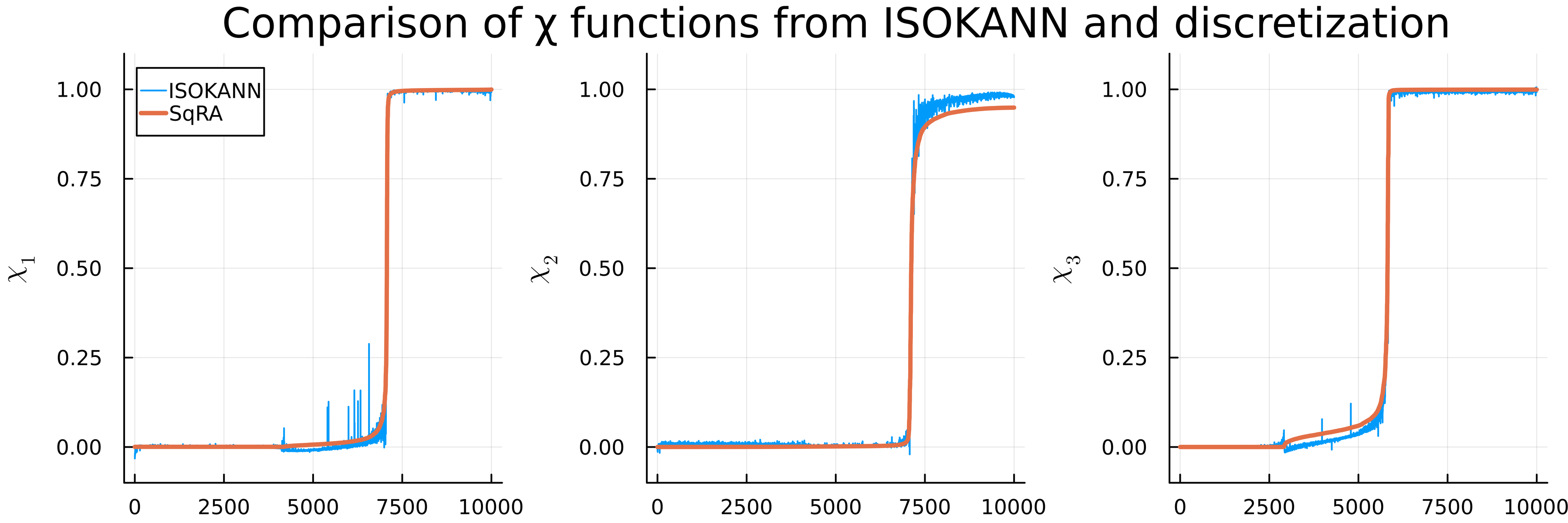}
    \caption{Comparison of the three components of $\chi: \mathbb{R}^3 \rightarrow \mathbb{R}^3$ obtained from either ISOKANN (blue) or direct discretization through the SqRA (red) for the potential from Section \ref{sec:ex3}.
    Each plot displays one component of $\chi$ evaluated at one of $\mathrm{10\,000}$ training samples, where the samples displayed on the horizontal axis were sorted by the respective $\chi$ value (different for each subfigure) obtained from SqRA to allow for a visual comparison. The random fluctuations that are seemingly visible for the $\chi$-functions computed via ISOKANN are the values on samples rather far away from all other samples; the gradients of the $\chi$-functions are still smooth functions and agree well for SqRA and ISOKANN.  
   The difference in scale of the second component results from details of the ISA algorithm and do not have a significant impact on the effective dynamics.
   }
    \label{fig:comparison}
\end{figure}

In these examples, we discretize the infinitesimal generator of the full dynamics using the Square-Root Approximation (SqRA) \cite{Lie2013, Donati2018b, Donati2021}, and then compute the $\chi$-functions using the PCCA+ algorithm \cite{Deuflhard2004, Kube2007, Weber2018}. 
The respective $\chi$-function computed via ISOKANN are $\chi$ functions conceptually equivalent and lead to essentially the same effective dynamics, see Fig.~\ref{fig:comparison}. However, the ISOKANN $\chi$-functions are only given on sample points and do not allow for convenient illustration. This approach isolates the errors introduced by the effective dynamics, providing a clear benchmark for our ISOKANN-based reduced models. For higher-dimensional systems, the $\chi$-functions are typically approximated from short burst simulations by neural networks via ISOKANN. %A comparison of the effective diffusion terms has shown that they agree well in the dynamically relevant regions.
% \AS{I tried to quantify this in a scalar like relative error, but could not figure out a single representative way.
% Most of the invariant measure accumulates in the "corner pixel" of the simplex (95\%). Here the diffusion is almost 0 such that the relative error is huge, even though the effect on the effective dynamics is negligible (here it is dominated by the drift).
% I think the effect is biggest in transition regions, where we indeed have the high diffusion, but the weight is very small.
% Indeed computing the "mean local relative error"
% \begin{align}
%     \int_\Delta \frac{ \|\tilde D(z)- D(z)\|}{\|D(z)\|} \mu (z) dz \approx 30
% \end{align}
% is huge. Going for the "global relative error" ($p=2$)
% \begin{align}
%     {\left(\frac{
%     \int_\Delta  \|\tilde D(z)- D(z)\|^p \mu (z) dz
%     }{
%     \int_\Delta \|D(z)\|^p \mu (z) dz
%     }\right)^{1/p}}  \approx 0.3
% \end{align}
% Not necessarily small.. But again, how meaningful is it?
% }

\subsection{One-dimensional system}
\begin{figure}[H]
    \centering
    \includegraphics[width=1\linewidth]{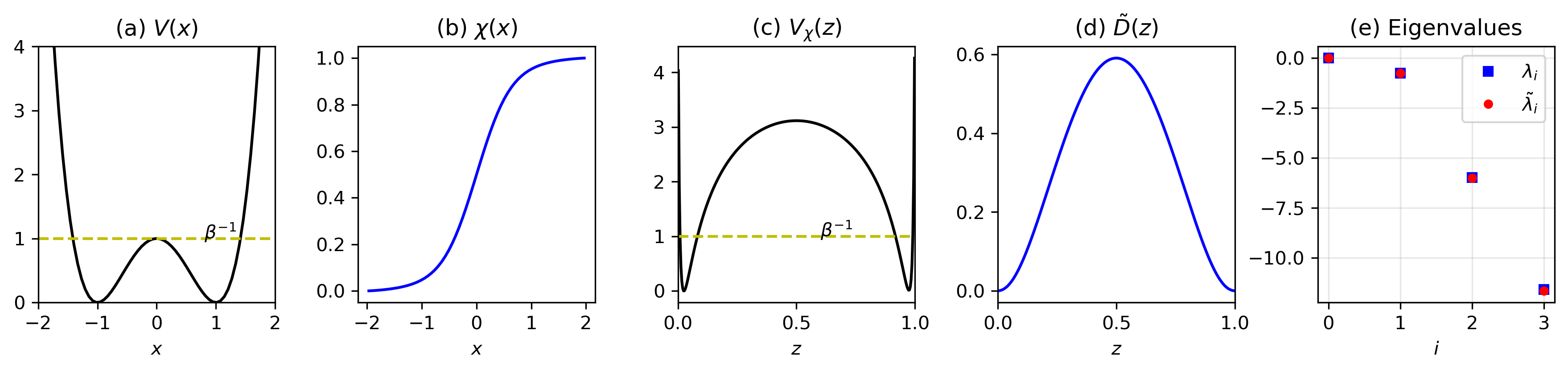}
        \vspace{-1cm}
    \caption{1D system. (a) Potential energy function; (b) $\chi$-function; (c) Effective potential; (d) Effective diffusion; (e) Eigenvalues.}
    \label{fig:fig3}
\end{figure}

As a first benchmark, we consider the overdamped Langevin dynamics in the one-dimensional ($d=1$) potential energy function
\[
V(x) = (x^2 - 1)^2, \quad \beta = 1,
\]
represented in Fig.~\ref{fig:fig3}, which provides an analytically tractable reference model.
The full generator $\cL$ was discretized on the uniform grid $[-2,2]$, partitioned into 63 equal intervals of width $\Delta x = 0.063$ using SqRA.
The associated $\chi$-function shown in Fig.~\ref{fig:fig3}(b) was obtained from PCCA+ assuming $m=1$.
The $\chi$-function smoothly interpolates between the two metastable basins, closely matching the true committor of the process. 

Using the $\chi$-function as collective variable, we computed the effective potential $V_{\chi}(z)$ via equation \eqref{eff-dyn-chi}. The result is displayed in Fig.~\ref{fig:fig3}(c), which shows that $V_\chi$ exhibits the double-well structure of the original potential $V(x)$ in latent space, but appears deformed due to the nonlinear transformation induced by the collective variable. 
While $V(x)$ describes the mechanistic energy landscape in full space, $V_{\chi}(z)$ represents the free-energy profile in latent space that incorporates both the enthalpic and entropic contributions of the system.

In Fig.~\ref{fig:fig3}(d), we show the effective diffusion coefficient $\tilde{D}(z)$ computed via \eqref{eq:DiffusionTensor}. 
It describes how thermal fluctuations are transmitted onto the latent space through the collective variable: lower values are observed within the metastable regions, where motion is constrained by the underlying potential, while higher values close to $z\approx 0.5$ indicate that barrier transitions are facilitated by increased entropic freedom.
Using $V_{\chi}(z)$ and $\tilde{D}(z)$, we discretized the effective generator defined in \eqref{eff-gen2} in latent space $[0,1]$ by applying SqRA.
The resulting eigenvalue spectrum of the original and the effective generators (Fig.~\ref{fig:fig3}(e)) shows agreement and exhibits a clear separation between the dominant slow mode and the remaining fast modes, consistent with the expected metastability of the system.
%

% %%%%%%%%%%%%%%% O L D

% We first consider the one-dimensional double-well potential 
% %
% \[
% V(x) = (x^2 - 1)^2, \quad \beta = 1,
% \]
% %
% with $\beta=1$ (Fig.~\ref{fig:fig3}) as an analytically tractable reference. 
% %
% We discretized the full generator $\cL$ (\ref{eq:full_generator}) on a uniform grid by SqRA [add references] and estimated the $\chi$-function using PCCA+ \cite{DeulfhardWeber}.
% %
% The $\chi$-function smoothly interpolates between the two metastable basins (Fig.~\ref{fig:fig3}(b)), closely matching the true committor of the process. 

% %
% The ISOKANN iteration recovers a $\chi$-function that smoothly interpolates between the two metastable basins (Fig.~\ref{fig:fig3}(b)), closely matching the full space committor of the process. 
% From this $\chi$, the effective potential $V_{\chi}(z)$ and diffusion coefficient $\tilde{D}(z)^2$ (Figs.~\ref{fig:fig3}(c) and (d)) can be computed directly via Eqs.~\eqref{eff-dyn-chi}--\eqref{eq:DiffusionTensor}. 
% The resulting eigenvalue spectrum of the original and the effective generators (Fig.~\ref{fig:fig3}(e)) shows agreements and exhibits a clear separation between the dominant slow mode and the remaining fast modes, consistent with the expected metastability.

\begin{figure}[H]
    \centering
    \includegraphics[width=1\linewidth]{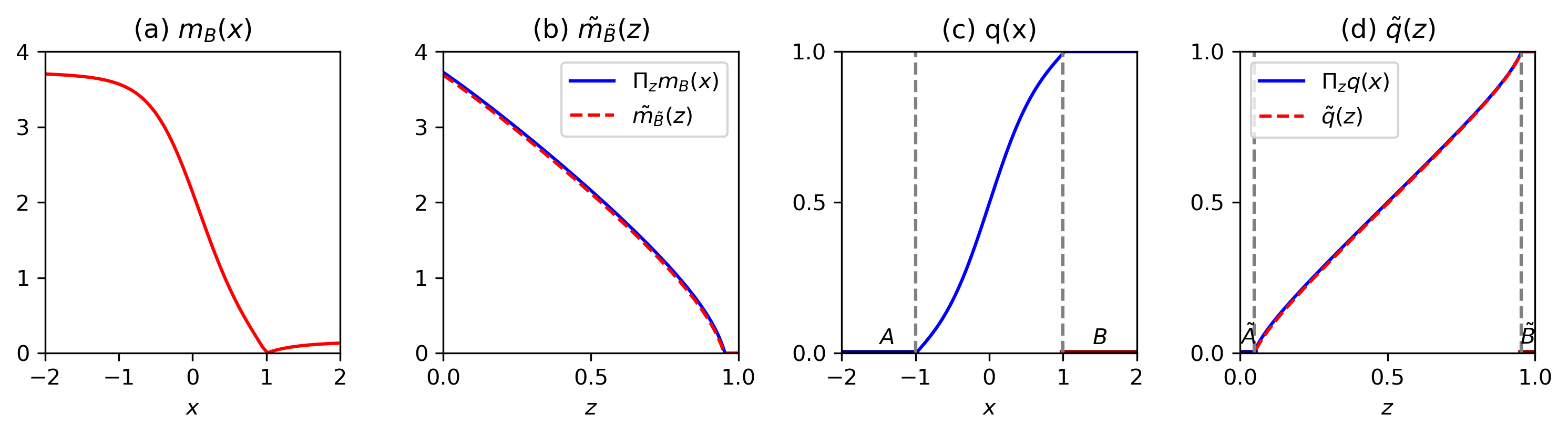}
        \vspace{-1cm}
    \caption{1D system. (a) MFPT in full space; (b) MFPT in latent space; (c) Committor in full space; (d) Committor in latent space;}
    \label{fig:fig4}
\end{figure}

To assess kinetic consistency, we compared MFPTs and committor functions computed for the full dynamics via \eqref{MFPT}--\eqref{committor} and the effective dynamics via \eqref{eff-MFPT-1d}--\eqref{eff-committor-1d}, respectively.
Fig.~\ref{fig:fig4}(a) shows the MFPT $m_B(x)$ in full space from a point $x$ to $x_B=1$, which corresponds to the right minimum of the potential, while Fig.~\ref{fig:fig4}(b) shows the projection of $m_B(x)$ onto the latent space and the effective MFPT $\tilde{m}_{\tilde{B}}$. 
Similarly, Figs.~\ref{fig:fig4}(c) and (d) show, respectively, the committor function in full space between the sets $A= \lbrace x\in \mathbb{R} : x\leq -1 \rbrace$ and $B = \lbrace x\in \mathbb{R} : x\geq 1 \rbrace$, and the corresponding projection in latent space together with the effective committor.
In both cases, the curves of the projected and effective quantities coincide, demonstrating exact agreement between the full dynamics and its effective representation.

The exact agreement between projected and full space quantities arises from the fact that, in one dimension ($d=1$), the projection onto the latent space reduces to a coordinate transformation, with the mapping $z=\chi(x)$ being invertible.
Consequently, no dynamical information is lost under this transformation, and the effective dynamics is fully equivalent to the original one.

In the following examples, we will show that such exact agreement is not always guaranteed.
The effective dynamics reproduces the kinetic behavior of the full system only when a clear spectral gap exists in the eigenvalue spectrum, ensuring a separation between slow and fast modes.

% Figure~\ref{fig:fig4} further compares the mean first-passage time (MFPT) and committor functions for the full (a) and effective dynamics (b).  Moreover, Fig.~\ref{fig:fig4} display the full committor $q$ (d) and the coomparison with the effective committor $\tilde{q}$ (e).
% The close agreement between the original and effective quantities confirms that ISOKANN accurately reproduces the effective kinetics in this minimal example.
%

\subsection{Two-dimensional system: enthalpic vs.\ entropic pathways}
\begin{figure}[H]
    \centering
    \includegraphics[width=1\linewidth]{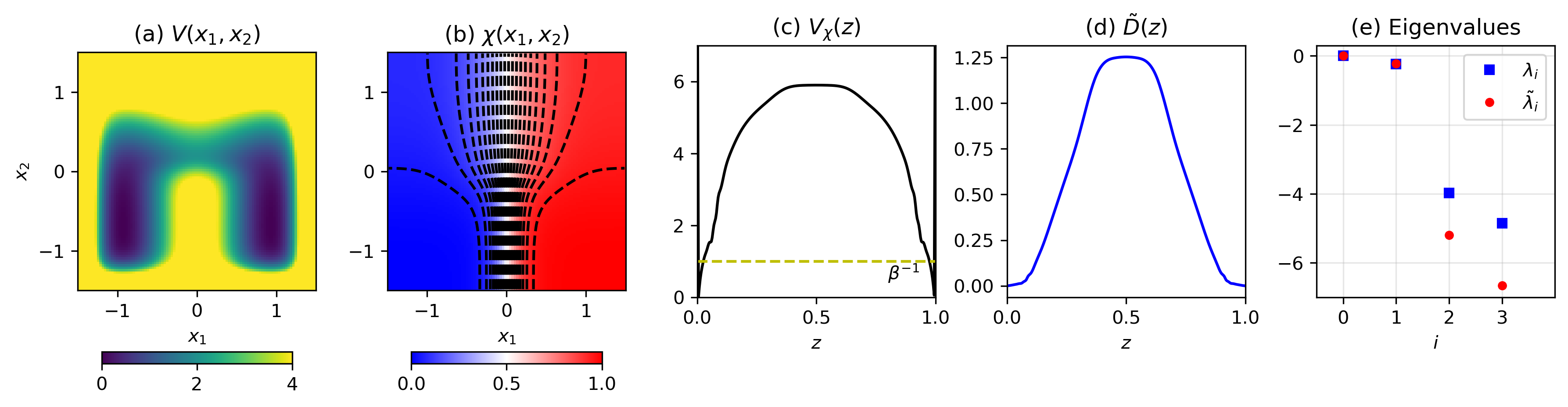}
        \vspace{-1cm}
    \caption{2D system. (a) Potential energy function; (b) $\chi$-function; (c) Effective potential; (d) Effective diffusion; (e) Eigenvalues.}
    \label{fig:fig5}
\end{figure}

As a second benchmark, we consider a two-dimensional ($d=2$) overdamped Langevin dynamics evolving on the potential energy surface
\begin{equation}
\begin{aligned}
V(x,y) = 0.2 \Big(
    &x^{12} + y^{12}
    + 20 \exp\big[-4x^{2} - 0.8\,(y+1)^{8}\big] \\[-3mm]
    &- 20 \ln\Big[
        \exp\big(-(x+1)^{2} - 0.1\,(y+1)^{4}\big)\\[-3mm]
    &    + \exp\big(-(x-1)^{2} - 0.1\,(y+1)^{4}\big)
    \Big]
\Big)\,,
\end{aligned}
\end{equation}
with $\beta=1$, visualized in Fig.~\ref{fig:fig5}(a).
The system exhibits two dominant metastable regions separated by an enthalpic barrier corresponding to a high-energy saddle in the potential energy function.
In addition, the two basins are connected by an entropic transition channel, i.e. by a region of the configuration space with high configuration multiplicity which raises the entropic contribution and consequently lowers the effective free energy barrier.
We discretized the full space $[-1.5,1.5]^2$ into 10.000 equal square bins of size $0.03\times 0.03$, then discretized the full generator $\cL$ using SqRA.
We estimated a single $\chi$-function, see Fig.~\ref{fig:fig5}(b), which identifies the dominant reaction coordinate between the metastable basins.

In Figs.~\ref{fig:fig5}(c-d) we show the effective potential $V_\chi(z)$ and the effective diffusion coefficient $\tilde{D}(z)$ built on the latent space $[0,1]$ defined by the $\chi$-function.
These quantities capture the energetic structure of the potential energy function and the interplay between the two dominant metastable sets.
However, since the $\chi$-function cannot distinguish whether transitions occur through the enthalpic barrier or via the entropic channel, the effective potential and diffusion coefficient should be interpreted as an averaged description of the dynamics, merging the contributions of both pathways into a single effective process.

Similar to the one-dimensional case, we discretized the effective generator $\tilde{\cL}$ by applying SqRA, using the effective potential and diffusion coefficient defined on the latent space.
The eigenvalue spectrum of the full and effective generators is shown in Fig.~\ref{fig:fig5}(e).
The two dominant eigenvalues are identical, indicating that the effective dynamics successfully reproduces the stationary distribution on the latent space and the principal transition between the metastable basins.
In contrast, the eigenvalues associated with faster modes show significant deviations, as expected from the dimensionality reduction (cf. \ref{order}).
%

%%%%%%%%%%%%%%%
% featuring full space dimension $d=2$, and build a one-dimensional latent space with $m=1$ (Fig.~\ref{fig:fig5}(a)). 
% %
% Here, the single $\chi$-function (Fig.~\ref{fig:fig5}(b)) resolves the dominant reaction channel, while the associated effective potential (Fig.~\ref{fig:fig5}(c)) and diffusion tensor (Fig.~\ref{fig:fig6}(d)) capture the interplay between the two dominant metastable sets. 
% The eigenvalue spectrum (Fig.~\ref{fig:fig5}(e)) shows that two dominant modes corresponding to transitions between metastable basins are (almost) identical for the full and effective generators, while the next eigenvalues (the ones representing faster modes) do not agree (as expected, cf.(\ref{order})).
%

\begin{figure}[H]
    \centering
    \includegraphics[width=1\linewidth]{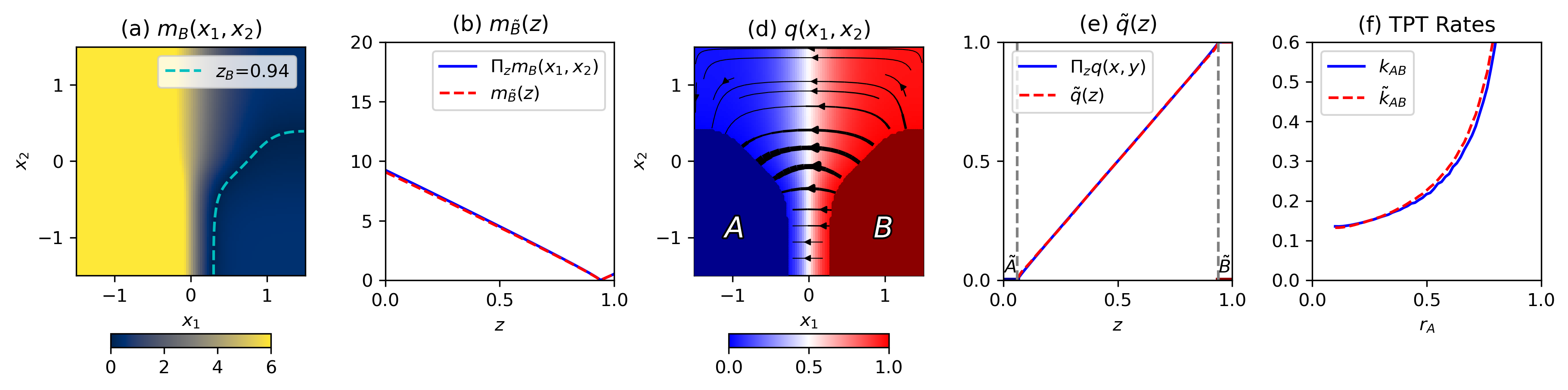}
        \vspace{-1cm}
    \caption{2D system. (a) MFPT in full space; (b) MFPT in latent space; (c) Committor in full space; (d) Committor in latent space; (f) TPT rates of full and effective dynamics.}
    \label{fig:fig6}
\end{figure}

We compared MFPTs and committor functions for the full and effective dynamics.
The MFPT in full space from a point $(x_1,x_2)$ to the set of points $B=\chi^{-1}(z_B=0.94)$ is displayed in Fig.~\ref{fig:fig6}(a), while its projection onto the latent space and the effective MFPT estimated by equation (\ref{eff-MFPT-1d}) is shown in Fig.~\ref{fig:fig6}(b).
Similarly, the committor functions are shown in in Figs.~\ref{fig:fig6}(c-d), with the sets
\[
\begin{aligned}
    \tilde{A} = \lbrace z \in[0,1]\, |\, z < 0.06 \rbrace, \quad\tilde{B} = \lbrace z \in[0,1]\, |\, z > 0.94 \rbrace,
\end{aligned}
\]
in latent space, and their corresponding pre-images $A=\chi^{-1}(\tilde{A})$, and $B=\chi^{-1}(\tilde{B})$ in full space.
This choice was made deliberately to investigate the effect of different transition pathways (enthalpic versus entropic) in situations where the latent space is unable to resolve their physical nature.
Despite this limitation, the effective dynamics reproduces the projected committor with high accuracy.

Finally, Fig.~\ref{fig:fig6}(f) compares the full TPT transition rate $k_{AB}$, as defined in (\ref{eq:full_rate}), to the effective rate $\tilde{k}_{AB}$ from (\ref{eff-rate}), as function of $r_A$, where $\tilde{A}=[0,r_A]$ and $\tilde{B}$ as above.
The close agreement between $k_{AB}$ and $\tilde{k}_{AB}$ confirms that the single collective variable $\chi$ retains the essential kinetic features of the two-dimensional process, even when distinct physical transition mechanisms coexist.

This example illustrates a scenario in which the projection is performed from a multidimensional system ($d>1$) onto a one-dimensional latent space ($m=1$).
Unlike the previous case, the transformation $z=\chi(x_1,x_2)$ is no longer invertible, and the projection entails a loss of information.
Nevertheless, when the eigenvalue spectrum exhibits a pronounced gap between the second and third eigenvalues $\lambda_1$ and $\lambda_2$, the effective dynamics can still reproduce the relevant kinetic quantities such as MFPTs and TPT rates between metastable regions. 
%

% Figures~\ref{fig:fig6} display the MFPT ((a) and (b)) and committor functions ((d) and (e)) in both full and latent space. The committor sets are defined as 
%
% \[
% \begin{aligned}
%     \tilde{A} = \lbrace z \in[0,1]\, |\, z < 0.06 \rbrace, \quad\tilde{B} = \lbrace z \in[0,1]\, |\, z > 0.94 \rbrace,
% \end{aligned}
% \]
%
% in latent space, and $A=\chi^{-1}(\tilde{A})$, and $B=\chi^{-1}(\tilde{B})$ in full space.
% This choice has been made deliberately in order to study the effect of different transition pathways (enthalpic vs.\ entropic) in cases where the latent space cannot resolve the difference. Furthermore, Fig.~\ref{fig:fig6}(f) compares the full TPT transition rate, $k_{AB}$, to the one resulting from the effective dynamics, $\tilde{k}_{AB}$, as function of $r_A$, with $\tilde{A}=[0,r_A]$ and $\tilde{B}$ as above.
% The effective dynamics accurately reproduces the mean first passage times as well as the TPT transition rates, confirming that the CV $\chi$ allows to retain the essential features of multidimensional transition pathways.
%

\subsection{Three-dimensional system with multiple $\chi$-functions}
\label{sec:ex3}
\begin{figure}[H]
    \centering
    \includegraphics[width=1\linewidth]{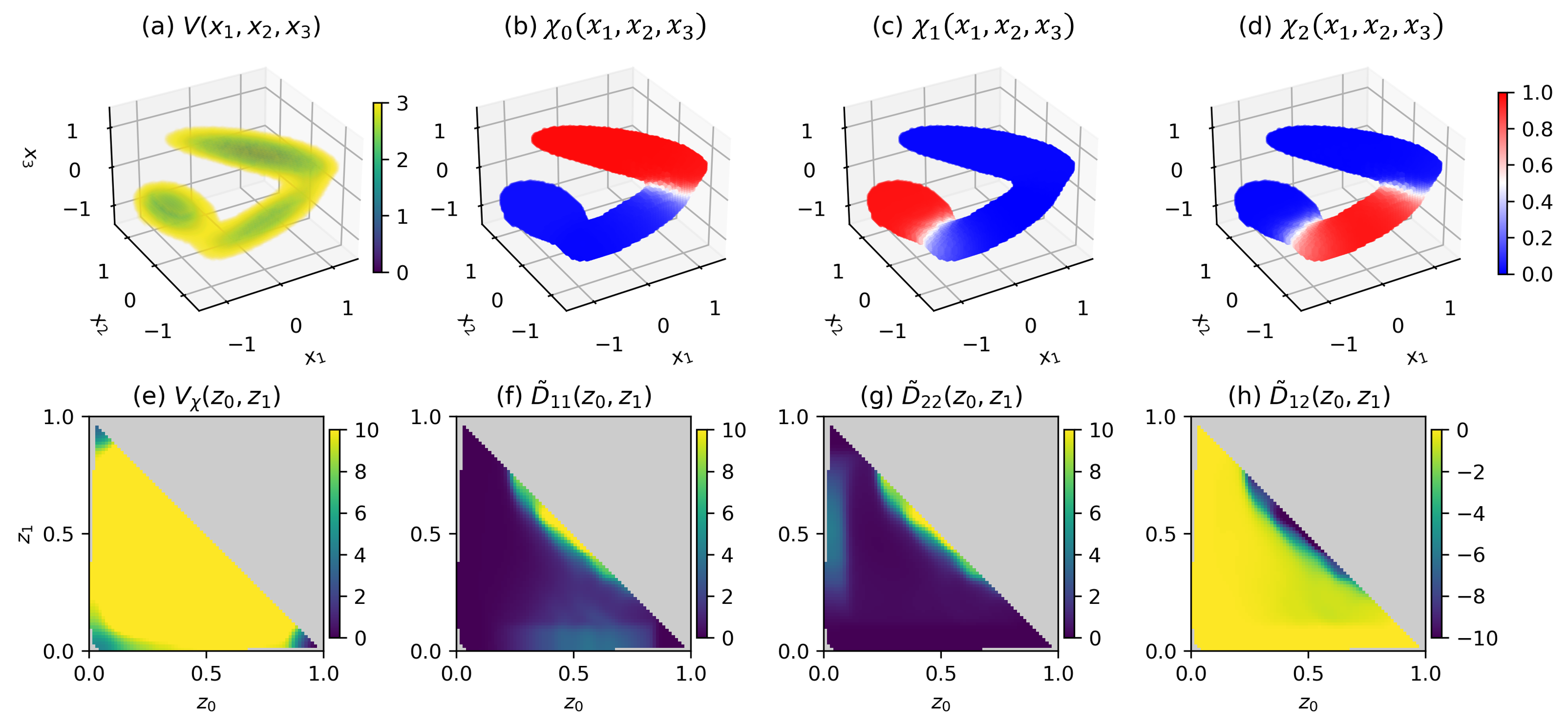}
        \vspace{-1cm}
    \caption{3D system. (a) Potential energy function; (b,c,d) Three $\chi$-functions; (e) Effective potential on $\bar{\bbDelta}^{2}$; (f-h) Diffusion tensor entries on $\bar{\bbDelta}^{2}$.}
    \label{fig:fig7}
\end{figure}

As final example, we consider a three-dimensional ($d=3$) overdamped Langevin dynamics with the potential energy function
\begin{align}
V(x,y,z) &= -15\,\exp\Big[
   -6.5\,(x-0.4)^2
   -11\,(x-0.4)(y-0)
   -5.5\,(y-0)^2
   -6.5\,(z-0.6)^2
\Big] \nonumber\\%[4pt]
&\quad -10\,\exp\Big[
   -1.5\,(x-0.2)^2
   -15\,(y+0.8)^2
   -5\,(z+0.4)^2
   +4\,(x-0.2)(z+0.4)
\Big] \nonumber\\%[4pt]
&\quad -12\,\exp\Big[
   -3\,(x+1)^2
   -1\,(y-0.5)^2
   -16.5\,(z+0.6)^2
\Big] \nonumber\\%[4pt]
&\quad +12\,\exp\Big[
   -1\,(x+0.5)^2
   -1\,(y-0.6)^2
   -10\,(z-0.1)^2
\Big] \nonumber\\%[4pt]
&\quad +0.8\,(x^4 + y^4 + z^4),
\end{align}
with $\beta=1$, shown in Fig.~\ref{fig:fig7}(a).
The system consists of three metastable regions, two of which are particularly deep and separated by a pronounced enthalpic barrier, while the third region acts as an intermediate metastable state that forms an entropic transition channel.
The full space $[-1.5,1.5]^3$ was discretized into 216.000 cubic bins of size $0.05 \times 0.05 \times 0.05$, and the generator $\cL$ was approximated using SqRA.
We estimated three $\chi$-functions via PCCA+ ($m=2$) shown in Figs.~\ref{fig:fig7}(b-d).
Removing the redundant coordinate $\chi_2=1-\chi_0-\chi_1$ yields the two-dimensional latent reduced simplex 
\[
    \bar{\bbDelta}^{2} = \left\lbrace z\in[0,1]^2 : z_0 + z_1 \le 1 \right\rbrace,
    \label{reduced-simplex2}
\]
on which the full dynamics is projected.
Using this two-dimensional latent space, we computed the effective potential $V_\chi(z_0, z_1)$ via (\ref{eff_pot_m_dim}), together with the components of the effective diffusion tensor $\tilde{D}_{11}(z_0,z_1)$, $\tilde{D}_{22}(z_0,z_1)$, and $\tilde{D}_{12}(z_0,z_1)$, as defined in (\ref{eq:DiffusionTensor}).
These quantities, reported in Figs.~\ref{fig:fig7}(e–h) respectively, reveal a structured free-energy landscape in the latent space $\bar{\bbDelta}^{2}$, with clear basins corresponding to the metastable sets of the original dynamics.
In contrast to the lower-dimensional cases, the diffusion tensor exhibits anisotropy and non-vanishing off-diagonal components, representing the coupling of thermal fluctuations along distinct transition channels.
%

% representing a system with three metastable states such that there exist $3$ dominant eigenvalues $\lambda_0, \lambda_1, \lambda_2$. In this case, we build a latent space with two dimensions ($m=2$), i.e., the effective description is built on two CVs $\chi_i \colon \X \to \R^2$, $i=1,2$, keeping in mind that the third $\chi$-function satisfies $\chi_3=1-\chi_1-\chi_2$, such that $\chi$-functions form a simplex partition of the state space; Figs.~\ref{fig:fig7} (b)-(d) show the resulting $\chi$-functions.  
% The resulting effective potential (Fig.~\ref{fig:fig7} (e)) and diffusion tensor entries (Fig.~\ref{fig:fig7} (f)-(h)) reveal a structured free-energy landscape in the latent space $(z_1, z_2)$, with clear basins corresponding to the metastable sets of the original dynamics.

\begin{figure}[H]
    \centering
    \includegraphics[width=1\linewidth]{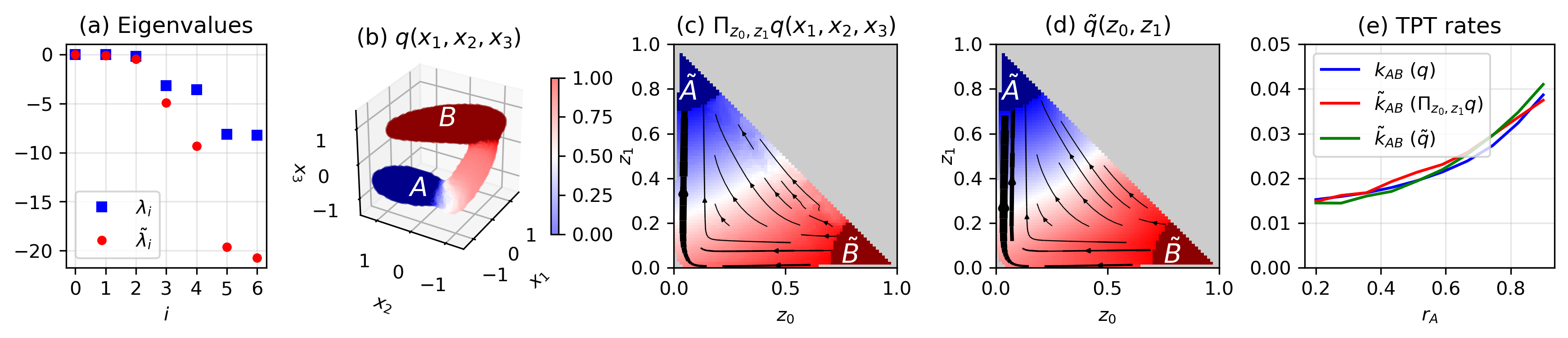}
    \vspace{-1cm}
    \caption{3D system. (a) Eigenvalues of $\cL$ and $\tilde{\cL}$; (b) Committor in full space; (c) Committor $q$ projected on $\bar{\bbDelta}^{2}$; (d) Effective committor $\tilde{q}$; (e) TPT rates of full and effective dynamics.}
    \label{fig:fig8}
\end{figure}

The spectrum of the effective generator $\tilde{\cL}$ shows three dominant eigenvalues that match those of the full dynamics, while the subsequent ones, $\tla_3,\tla_4,\ldots$, are consistently smaller than those of the full dynamics, as expected.  
The full-space committor function $q(x_1,x_2,x_3)$, its projection onto $\bar{\bbDelta}^{2}$, and the corresponding effective committor $\tilde{q}(z_0,z_1)$, evaluated by solving \eqref{committor_m}, are displayed in Figs.~\ref{fig:fig8}(b-d), 
where the committor sets in $\bar{\bbDelta}^{2}$ are defined as 
\[
\begin{aligned}
\tilde{A} &= \bigl\{ (z_0,z_1)\in[0,1]^2 \,\big|\, \lVert (z_0,z_1)-(0,1)\rVert_2 \le 0.3 \bigr\},\\
\tilde{B} &= \bigl\{ (z_0,z_1)\in[0,1]^2 \,\big|\, \lVert (z_0,z_1)-(1,0)\rVert_2 \le 0.3 \bigr\},
\end{aligned}
\]
with $A=\chi^{-1}(\tilde{A})$ and $B=\chi^{-1}(\tilde{B})$ corresponding pre-images in full space.
Figs.~\ref{fig:fig8}(c) and \ref{fig:fig8}(d) also display the reactive fluxes, $j_{AB}$ and $\tilde{j}_{AB}$, highlighting the main reaction channels.
Comparing Figs.~\ref{fig:fig8}(c) and \ref{fig:fig8}(d) demonstrates that the coarse-grained model retains both the topology and the transition ordering of the underlying process. 

Finally, Fig.~\ref{fig:fig8}(e) compares the full TPT transition rate $k_{AB}$ with the effective rate $\tilde{k}_{AB}$ as a function of the boundary parameter $r_A$ in the definition of the set $\tilde{A}$.
The close agreement between the two rates confirms that the reduced dynamics not only reproduces the correct spatial structure of the transition but also preserves the associated kinetics.
Small discrepancies between full TPT rates $k_{AB}$ and effective TPT rates $\tilde{k}_{AB}$ arise from the combined effect of dimensionality reduction, which introduces memory effects neglected in the effective Markovian representation of the system, and numerical approximations induced by the discretization of space and the calculation of the effective potential and the diffusion tensor components.
The fact that these differences remain modest indicates that the reduced model successfully captures the dominant kinetic behavior. 
%
%

%These results illustrate that ISOKANN systematically identifies low-dimensional latent variables which encode metastable connectivity and yield quantitatively accurate effective kinetics.

% \begin{figure}[H]
%     \centering
%     \includegraphics[width=1\linewidth]{fig8.png}
%     \caption{3D system. (a) Eigenvalues; (b) Committor $q$; (c) Committor $q$ projected on $(z_1, z_2)$; (d) Committor $\tilde{q}$.}
%     \label{fig:fig8}
% \end{figure}

\subsection{Discussion}

Across all examples, the CV $\chi$ of ISOKANN allows to consistently reconstruct the dominant eigenmodes and corresponding kinetic observables of the underlying stochastic process. 
The one-dimensional case validates the formal connection between $\chi$-functions, committor functions, and effective potential and diffusion tensor. 
The two- and three-dimensional systems highlight how the algorithm balances energetic and entropic contributions and generalizes to higher-dimensional latent spaces. 
Together, these numerical experiments confirm that Koopman-inspired learning of invariant subspaces provides a robust route to constructing effective stochastic models with faithful transition rates and pathways.

\section{Conclusions and Outlook}\label{sec:conclusion}

This work presents a unified framework that bridges data-driven learning of collective variables (CVs) with the theoretical foundations of effective stochastic dynamics. Building on transfer-operator and Koopman operator theory, the authors introduce the ISOKANN method as a systematic approach to extract low-dimensional variables that preserve the slow dynamical modes of complex stochastic systems, such as those encountered in molecular dynamics.

At its core, ISOKANN identifies nonlinear mappings, called $\chi$- or membership functions, that span the invariant subspace of the dominant eigenmodes of the Koopman operator. By combining neural-network representations with iterative power-like updates and geometric normalization (via the Inner Simplex Algorithm), the method discovers collective variables without prior physical intuition. The learned $\chi$-functions define a latent space in which an \emph{effective dynamics}, a projected Ornstein–Uhlenbeck process with state-dependent noise, can be derived analytically. This effective dynamics is Markovian and is guaranteed to inherit the dominant eigenvalues and eigenfunctions of the full generator, ensuring that metastable transitions and relaxation timescales are faithfully reproduced.

The study further links this operator-theoretic perspective to transition path theory (TPT), showing how committor functions, mean first-passage times, reactive fluxes, and transition rates can be computed consistently in both the original and effective dynamics. This connection provides a rigorous kinetic interpretation of the ISOKANN-derived effective dynamics.

Through numerical experiments on benchmark one-, two-, and three-dimensional potential landscapes, the authors demonstrate that effective dynamics based on ISOKANN robustly reconstructs the dominant reaction coordinates and reproduces transition rates and pathways across both enthalpic and entropic barriers. The approach accurately captures metastable connectivity and reactive flux structure while reducing the dimensionality of the problem.

Future research on ISOKANN will have to tackle two main challenges:  its extension to higher-dimensional molecular systems and to non-reversible processes. While the theoretical formulation of ISOKANN directly generalizes to high-dimensional input spaces, practical applications reveal that the main obstacle is the speed of convergence of the iterative ISOKANN algorithm when identifying collective variables in very large state spaces. Improving convergence efficiency, through accelerated subspace updates and adaptive sampling, is a current focus of our ongoing research.
With regard to non-reversible dynamics, such as Langevin systems in joint position–velocity space, recent theoretical advances already provide a pathway: \cite{ZhangSchuette2023Understanding} includes first rigorous results describing how effective dynamics and spectral subspaces can be consistently defined for such cases. These results lay the foundation for extending ISOKANN to general molecular processes beyond the reversible setting considered here.

A further direction of development concerns adaptive collective-variable learning and feedback coupling. The general aim is to iteratively refine the present CVs by sampling along the transition pathways generated by the effective dynamics they induce. This feedback loop couples the learning of $\chi$-functions with guided sampling in latent space: new trajectory data are drawn preferentially from regions of high reactive flux or uncertainty, and the resulting updated ISOKANN model improves the resolution of transition channels and rare events. Such adaptive coupling between model learning and dynamical exploration promises a self-consistent route toward efficient discovery of collective variables and accurate reconstruction of effective kinetics in complex, high-dimensional systems.

In summary, ISOKANN offers a principled and interpretable machine-learning framework for stochastic model reduction: it merges neural representation learning with Koopman spectral theory and stochastic reduction analysis. The results underline that algorithmic discovery of CVs and the consistent formulation of the induced effective dynamics can advance the quantitative understanding of complex molecular systems, opening perspectives toward adaptive learning, non-reversible dynamics, and data-driven coarse-grained modeling in high dimensions.

\subsection*{Acknowledgments}
This work has been partially funded by the Deutsche Forschungsgemeinschaft (DFG, German Research Foundation) through grant CRC 1114/2-3 (Projects B03, A04, A05, and B05) and under Germany’s Excellence Strategy -- The Berlin Mathematics Research Center MATH+ (EXC-2046/1 project ID:  390685689).

%%%%%%%%%%%%%%%%%%%%%%%%%%%%%%%%%%%%%%%%%%%%%%%%%%%%%%%%%%
%%%%%%%%%%%%%%%%%%%%%%%%%%%%%%%%%%%%%%%%%%%%%%%%%%%%%%%%%%
%%%%%%%%%%%%%%%%%%%%%%%%%%%%%%%%%%%%%%%%%%%%%%%%%%%%%%%%%%
%%%%%%%%%%%%%%%%%%%%%%%%%%%%%%%%%%%%%%%%%%%%%%%%%%%%%%%%%%
%%%%%%%%%%%%%%%%%%%%%%%%%%%%%%%%%%%%%%%%%%%%%%%%%%%%%%%%%%
\appendix
\section{Derivation of the effective potential}
In what follows, we derive the effective potential for the effective dynamics.

In the case $m=1$, the effective Fokker-Planck operator acting on probability density functions $\trho(z)$ defined on the latent space $\mathbb{Z}=[0,1]$ reads
\begin{equation}
\begin{aligned}
\tLast \trho(z) 
&=
- \frac{\partial\big((a + \lambda z)\,\trho(z)\big)}{\partial z} + \frac{1}{2} \frac{\partial^2 \big(\hsigma^2(z)\,\trho(z)\big)}{\partial z^2} \\ \label{FP_current}
&= - \frac{\partial \tilde{J}(z)}{\partial z},
\end{aligned}
\end{equation}
where we introduced the \emph{probability current}
\[
\tilde{J}(z) = (a + \lambda z)\,\trho(z)  - \frac{1}{2} \frac{\partial \big(\hsigma^2(z)\,\trho(z)\big)}{\partial z}.
\]
Stationary solutions of (\ref{FP_current}) satisfy $\frac{\partial \tilde{J}}{\partial z} = 0$, implying that the current $\tilde{J}(z)$ is constant on $[0,1]$. 
Imposing reflecting boundary conditions at $z=0$ and $z=1$ gives 
$\tilde{J}(0)=\tilde{J}(1)=0$, and thus $\tilde J(z)= 0$ on $[0,1]$.
Consequently, the unique invariant density of the effective dynamics is given by
\[
\tilde{\mu}(z)=\frac{1}{Z_\chi} \exp(- V_\chi(z)),\qquad Z_\chi=\int_0^1 \exp(- V_\chi(z))\,dz,
\]
where we introduced the effective potential
\begin{eqnarray}
V_{\chi}(z) & = & \log\Big(\frac{1}{2}\hsigma(z)^2\Big) - \phi(z)\,+\,\text{const},
\\
\phi(z) & = & \int_{z^\ast}^{z} 2 \frac{a + \lambda y}{\hsigma(y)^2} \, dy .
\end{eqnarray}

In the case $m>1$, the effective Fokker-Planck operator on  $\bar{\bbDelta}^{m}$ takes the form
\begin{equation}
\label{eff_FP_m_dim}
\begin{aligned}
\tLast \trho(z) 
&= 
- \nabla_z \cdot \bigl((b + Qz)\,\tilde\rho(z)\bigr)
+ \frac{1}{2}
\sum_{k,\ell=1}^m
\frac{\partial^2}{\partial z_k \partial z_\ell}
\Bigl( \bigl(\hsigma\hsigma^\top\bigr)_{k\ell}(z)\,\tilde\rho(z) \Bigr),
\end{aligned}
\end{equation}
where we introduced the probability current
\[
\tilde{J}(z) =
(b+Qz)\,\trho(z)  - 
\frac{1}{2} 
\nabla_{z} 
\cdot 
\left(
(\hsigma\hsigma^\top)(z)\,\trho(z)
\right).
\]
Unlike the $m=1$ case, stationarity with reflecting boundaries does not imply $\tilde{J}(z)=0$ on $\bar{\bbDelta}^{m}$.
To eliminate probability currents, we need reversibility, which is guaranteed by the self-adjointness of $\tL$ on $\tilde{\mathcal{H}}$.
Under this condition, $\tilde{J}(z)=0$ on  $\bar{\bbDelta}^{m}$, leading to
\begin{equation}
\label{eff_free_flux_m_dim}
(b+Qz)\,\trho(z) 
= 
\frac{1}{2} 
\nabla_{z} 
\cdot 
\left(
(\hsigma\hsigma^\top)(z)\,\trho(z)
\right),
\end{equation}

Under dimension reduction ($m+1\rightarrow m$), the matrix $\hsigma\hsigma^{\top}$ defined in (\ref{eq:DiffusionTensor}) is positive definite, and therefore invertible. 
Thus, (\ref{eff_free_flux_m_dim}) can be rewritten as
\begin{equation}
\begin{aligned}
-
\nabla_z \log \tmu(z)
&=
2 (\hsigma\hsigma^\top)^{-1} 
\left(
\frac{1}{2}
\nabla_{z} \cdot (\hsigma\hsigma^\top)
-
(b+Qz)(z)
\right) \\
% &= 
% 2 (\hsigma\hsigma^\top)^{-1} 
% (b+Qz)_l(z)
% -
% (\hsigma\hsigma^\top)^{-1} 
% \nabla_z \cdot (\hsigma\hsigma^\top)
% \\
&=
\psi(z).
\end{aligned}
\end{equation}
%
% Componentwise alternative:
% %
% \begin{equation}
% \begin{aligned}
% \frac{\partial \log \tmu}{\partial z_l}
% &=
% 2 
% \sum_{j=0}^{m-1}
% \left((\hsigma\hsigma^\top)^{-1} \right)_{lj}
% \left(
% (b+Qz)_l(z)
% -
% \frac{1}{2}
% \sum_{k=0}^{m-1} 
% \frac{\partial (\hsigma\hsigma^\top)_{jk}}{\partial z_k}
% \right) 
% \\
% &=g_l(z).
% \end{aligned}
% \end{equation}
%
Since the effective dynamics is reversible, the vector field $\psi(z)$ is curl–free on $\bar{\bbDelta}^{m}$:
\[
\frac{\partial \psi_k}{\partial z_l}=\frac{\partial \psi_l}{\partial z_k}.
\]
Hence, there exists an effective invariant density
\[
\tilde{\mu}(z)=\frac{1}{Z_\chi} \exp(- V_\chi(z)),\qquad Z_\chi=\int_{\bar{\bbDelta}^{m}} \exp(- V_\chi(z)) \, dz,
\]
where the effective potential is defined, up to an additive constant, by
\begin{equation}
\begin{aligned}
V_{\chi}(z) = 
\int_{\gamma:z^\ast\to z}
\psi(y) \cdot dy,
\end{aligned}
\end{equation}
with the line integral taken along any smooth path $\gamma\subset \bar{\bbDelta}^{m}$ from a starting point $z_\ast$ to $z$. 
The curl-free condition ensures that the effective potential is path-independent.

A practical strategy is to evaluate the integral from $z^\ast$ to $z$ coordinate by coordinate along the axes:
\[
V_\chi(z)
=
V_\chi(z_\ast)
+
\sum_{l=0}^{m-1}
\int_{z_l^\ast}^{z_l}
\psi_i\!\left(z_1',\dots,z_{l-1},\,
y_l,
\,z_{l+1}^\ast,\dots,z_m^\ast
\right)
\,dy_l.
\]
For example, for the case $m=2$, we obtain
\[
    V_{\chi}(z_0,z_1)
    =
    V_\chi(z_0^\ast, z_1^\ast)
    +
    \left(
    \int_{z_0^\ast}^{z_1} \psi_0(y_0,z_1^\ast) \, dy_0
    +
    \int_{z_1^\ast}^{z_1} \psi_1(z_0^\ast,y_1) \, dy_1
    \right)
\]
Since $V_\chi$ is defined up to an additive constant, we can set $V_\chi(z_0^\ast, z_1^\ast)=0$.
\bibliographystyle{unsrt}
\bibliography{refs}

\end{document}